\documentclass[11pt]{article}
\usepackage[pdftex]{graphicx}
\usepackage[inline]{enumitem}
\usepackage{mathtools}
\usepackage{amssymb}
\usepackage{amsthm}
\usepackage{dsfont}
\usepackage{chngcntr}
\usepackage{multirow}
\usepackage{graphicx}
\usepackage{caption}
\usepackage{subcaption}
\usepackage{longtable}
\usepackage{xcolor}

\newtheorem{thm}{Theorem}[section]
\newtheorem{lemma}{Lemma}[section]
\newtheorem{corollary}{Corollary}[section]
\theoremstyle{remark}
\newtheorem{ass}{Assumption}
\theoremstyle{definition}
\newtheorem{remark}{Remark}[section]

\newcounter{step}
\newenvironment{steps}
{\begin{list}{\it{Step\,\arabic{step}} :}
{\usecounter{step}}}
{\end{list}}

\newcommand*{\rom}[1]{\expandafter{\romannumeral #1\relax}}

\newcommand{\E}{\mathbb{E}}
\newcommand{\Var}{\mathrm{Var}}

\allowdisplaybreaks

\renewcommand{\baselinestretch} {1.2}

\makeatletter 
\def\singlespace{\def\baselinestretch{1}\@normalsize}

\title{{\sc Moving Block and Tapered Block Bootstrap for Functional Time Series with an Application to the $K$-Sample Mean Problem}}

\author{ {Dimitrios ~{\sc PILAVAKIS}}, \; {Efstathios ~{\sc PAPARODITIS}\footnote{
          Corresponding author (email: \texttt{stathisp@ucy.ac.cy})}} \; and \;  {Theofanis ~{\sc SAPATINAS}} \\
          Department of Mathematics and Statistics, University of Cyprus, \\
          P.O. Box 20537, CY 1678 Nicosia, CYPRUS.
}
\date{}
\begin{document}
\maketitle

\begin{abstract}
We consider infinite-dimensional Hilbert space-valued random variables that are assumed to be temporal dependent in a broad sense. We prove a central limit theorem for the moving block bootstrap and for the tapered block bootstrap, and show that these block bootstrap procedures also provide consistent estimators of the long run covariance operator. 
Furthermore, we consider block bootstrap-based procedures for fully functional testing of the equality of mean functions between several independent functional time series. 
We establish validity of the block bootstrap methods in approximating the distribution of the statistic of interest under the null and show consistency of the block bootstrap-based tests under the alternative. The finite sample behaviour of the procedures is investigated by means of simulations. 
An application to a real-life dataset is also discussed.

\medskip
\noindent
{\em Some key words:} {\sc
 Functional Time Series; Mean Function; Moving Block Bootstrap; Tapered Block Bootstrap; Spectral Density Operator; $K$-sample mean problem}

\end{abstract}

\section{\sc Introduction}

In statistical analysis, conclusions are commonly derived based on information obtained from a random sample of observations. In an increasing number of fields, these observations are curves or images which are viewed as functions in appropriate spaces, since an observed intensity is available at each point on a line segment, a portion of a plane or a volume. Such observed curves or images are called `functional data'; see, e.g., Ramsay and Dalzell (1991), who also introduced the term `functional data analysis' (FDA) which refers to statistical methods used for analysing this kind of data.
\par
In this paper, we consider observations stemming from a stochastic process $\mathds{X}=(X_t,\,t\in \mathbb{Z})$ of Hilbert space-valued random variables which satisfies certain stationarity and weak dependence properties. Our goal is to infer properties of the stochastic process based on an observed stretch $X_1,X_2,\ldots,X_n,$ i.e., on a functional time series. In this context, we commonly need to calculate the distribution, or parameters related to the distribution, of some statistics of interest based on $X_1,X_2,\ldots,X_n$. Since in a functional set-up such quantities typically depend in a complicated way on infinite-dimensional characteristics of the underlying stochastic process $\mathds{X}$, their calculation is difficult in practice. As a result, resampling methods and, in particular, bootstrap methodologies are very useful.
\par
For the case of independent  and identically distributed (i.i.d.) Banach space-valued random variables, Gin\'{e} and Zinn $(1990)$ proved the consistency of the standard i.i.d. bootstrap for the sample mean. For functional time series, Politis and Romano $(1994)$ established validity of the stationary bootstrap for the sample mean and for (bounded)  Hilbert space-valued random variables satisfying certain mixing conditions. A functional sieve bootstrap procedure  for functional time series has been proposed by Paparoditis (2017). Consistency of the non-overlapping block bootstrap for the sample mean and for near epoch dependent Hilbert space-valued random variables has been established by Dehling {\em et al.} (2015). However, up to date, consistency results are not available for the moving block bootstrap (MBB) or its improved versions, like the tapered block bootstrap (TBB), for functional time series.
Notice that the MBB for real-valued time series was introduced by K\"{u}nsch $(1989)$ and Liu and Singh $(1992)$. The basic idea is to resample blocks of the time series and to joint them together in the order selected in order  to form a new set of pseudo observations. This resampling scheme retains the dependence structure of the time series within each block and can be, therefore, used to approximate  the distribution of a wide range of statistics. 
The TBB for real-valued time series was introduced  by Paparoditis and Politis (2001). It  uses a taper window to downweight the observations  at the beginning and at the end of each resampled block and  improves the bias properties of the MBB.
\medskip

\par
The aim of this paper is twofold. First, we prove consistency of the MBB and of the TBB for the sample mean function in the case of weakly dependent Hilbert space-valued random variables. Furthermore, we show that these bootstrap methods provide consistent estimators of the covariance operator of the sample mean function estimator, that is of the spectral density operator of the underlying stochastic  process at frequency zero. We derive our theoretical results under quite general dependence assumptions on $\mathds{X},$ i.e., under $L^2$-$m$-approximability assumptions, which are satisfied by a large class of commonly used functional time series models; see, e.g., H\"{o}rmann and  Kokoszka (2010). Second, we apply the above mentioned bootstrap procedures to the problem of fully functional testing of the equality of the mean functions between a number of independent functional time series.
Testing the equality of mean functions for i.i.d. functional data has been extensively discussed in the literature; see, e.g., Benko {\em et al.} (2009), H\'{o}rvath and Kokoszka (2012, Chapter 5), Zhang (2013) and Staicu {\em et al.} (2015). 
Bootstrap alternatives over asymptotic approximations have been proposed in the same context by Benko {\em et al.} (2009), Zhang {\em et al.} (2010) and, more recently, by Paparoditis and Sapatinas (2016).
Testing equality of mean functions for dependent functional data has also attracted some interest in the literature. Horv\'{a}th {\em et al.} (2013) developed an asymptotic procedure for testing equality of two mean functions for functional time series. Since the limiting null distribution of a fully functional, $L^2$-type test statistic, depends on difficult to estimate process characteristics, tests are considered which are based on a finite number of projections.  A projection-based test has also been considered by Horv\'{a}th and Rice (2015). Although such tests lead to manageable limiting distributions, they have non-trivial power only for deviations from the null which are not orthogonal to the subspace generated by the particular projections considered.
\par
In this paper, we show that the MBB and TBB procedures can be successfully applied to approximate the distribution under the null of such fully functional test statistics. This is achieved by designing the suggested block bootstrap procedures in such a way that the generated pseudo-observations satisfy the null hypothesis of interest. Notice that such block bootstrap-based testing methodologies are applicable to a broad range of possible test statistics. As an example, we prove validity for the $L^2$-type test statistic recently proposed by Horv\'{a}th {\em et al.} (2013). 
\par
The paper is organised as follows. In Section~\ref{sec:BB}, the basic assumptions on the underlying stochastic process $\mathds{X}$ are stated and the MBB and TBB procedures for weakly dependent, Hilbert space-valued random variables, are described. Asymptotic validity of the block bootstrap procedures for estimating the distribution of the sample mean function is established and consistency of the long run covariance operator, i.e., of  the spectral density operator of the underlying stochastic process at frequency zero, is proven. Section~\ref{sec:testing} is devoted to the problem of testing equality of mean functions for several independent functional time series. Theoretical justifications of an appropriately modified version of the MBB and of the TBB procedure for approximating the null distribution of a fully functional test statistic is given and consistency under the alternative is shown. Numerical simulations  and a real-life data example are presented and discussed in Section~$4$. Auxiliary results and proofs of the main results are deferred to Section~$5$ and to the supplementary material.

\section{\sc Block Bootstrap Procedures for Functional Time Series}
\label{sec:BB}
\subsection{\sc Preliminaries and Assumptions}
\label{ssec:Preliminaries}
We consider a strictly stationary stochastic process $\mathds{X}=\{X_t,\,t\in \mathbb{Z}\},$ where the random variables $X_t$ are random functions $ X_t(\omega, \tau),\,\tau\in \mathcal{I},\,\omega\in\Omega,\,t\in\mathbb{Z}$, defined on a probability space $(\Omega,A, P)$ and take values in the separable Hilbert-space of squared-integrable $\mathbb{R}$-valued functions on $\mathcal{I}$, denoted by $L^2(\mathcal{I}).$ The expectation function of $X_t$, $\E X_t\in L^2(\mathcal{I}),$ is independent of $t,$ and it is denoted by $\mu.$ Throughout Section \ref{sec:BB}, we assume for simplicity that $\mu=0.$
We define $\langle f, g\rangle =\int_{\mathcal{I}}f(\tau)g(\tau)\mathrm{d}\tau,$ $\|f\|^2=\langle f,f\rangle$ and the tensor product between $f$ and $g$ by $f\otimes g(\cdot)=\langle f, \cdot\rangle g.$ For two Hilbert-Schmidt operators $\Psi_1$ and $\Psi_2$, we denote by $\langle \Psi_1,\Psi_2\rangle_{HS} =\sum_{i=1}^{\infty}\langle\Psi_1(e_i),\Psi_2(e_i)\rangle$ the inner product which generates the Hilbert-Schmidt norm $\|\Psi_1\|_{HS}=\sum_{i=1}^{\infty}\|\Psi_1(e_i)\|^2$, for $\{e_i,i=1,2,\ldots\}$ an orthonormal basis of $L^2(\mathcal{I}).$
Without loss of generality, we assume that $\mathcal{I} = [0, 1]$ (the unit interval) and, for simplicity, integral signs without the limits of integration imply integration over the interval $\mathcal{I}.$ We finally write $L^2$ instead of $L^2(\mathcal{I})$.
\par
To describe the dependent structure of the stochastic process $\mathds{X}$, we use the notion of $L^p$-$m$-approximability; see H\"{o}rmann and  Kokoszka (2010). A stochastic process $\mathds{X}=\{X_t,t\in\mathbb{Z}\}$ with $X_t$ taking values in $L^2$, is called $L^2$-$m$-approximable if the following conditions are satisfied:
\begin{enumerate}[label=(\roman*)]
     \item $X_t$ admits the representation \begin{equation}\label{orismosL2a}X_t=f(\delta_t,\delta_{t-1},\delta_{t-2},\ldots)\end{equation}  for some measurable function $f:S^\infty\rightarrow L^2$, where $\{\delta_t,\,t\in\mathbb{Z}\}$ is a sequence of i.i.d. elements in $L^2$.
\item $\E\|X_0\|^2<\infty$ and
    \begin{equation}\label{orismosL2b}
    \sum_{m\geq1}\sqrt{\E\|X_t-X_{t,m}\|^2}<\infty,\end{equation} where $X_{t,m}=f(\delta_t,\delta_{t-1},\ldots,\delta_{t-m+1},\delta_{t,t-m}^{(m)},\delta_{t,t-m-1}^{(m)},\ldots)$
    and, for each $t$ and $k,$ $\delta_{t,k}^{(m)}$ is an independent copy of $\delta_t.$
\end{enumerate}

The intuition behind the above definition is that the function $f$ in \eqref{orismosL2a} should be such that the effect of the innovations $\delta_i$ far back in the past becomes negligible, that is, these innovations can be replaced by other, independent, innovations. We somehow strengthen \eqref{orismosL2b} to the following assumption.
\begin{ass} \label{as:L2}
$\mathds{X}$ is $L^2$-$m$-approximable and satisfies
$$\lim_{m\to\infty}m\sqrt{\E\|X_t-X_{t,m}\|^2}=0.$$
\end{ass}
\noindent Notice that the above assumption is satisfied by many linear and non-linear functional time series models cconsidered in the literature; see, e.g.,  H\"{o}rmann and  Kokoszka (2010).
\medskip


\subsection{\sc The Moving Block Bootstrap}
\label{ssec:MBB}
The main idea of the MBB is to split the data into overlapping blocks of length $b$ and to obtain the bootstrapped pseudo-time series by joining together the $k$ independently and randomly selected blocks of observations in the order selected. Here, $k$ is a positive integer satisfying $b(k-1)<n$ and $bk\geq n.$ For simplicity of notation, we assume throughout the paper that $n=kb$. Since the dependence of the original time series is maintained within each block, it is expected that for weakly dependent time series, this bootstrap procedure will, asymptotically, correctly imitate the entire dependence structure of the underlying stochastic process if the block length $b$ increases to infinity, at some appropriate rate, as the sample size $n$ increases to infinity. Adapting this resampling idea to a functional time series $\boldsymbol{\mathrm{X}_n}=\{X_t,\,t=1,2,\ldots,n\}$ stemming from a strictly stationary stochastic process $\mathds{X}=\{X_t,t\in\mathbb{Z}\}$ with $X_t$ taking values in $L^2$ and $\E(X_t)=0$, leads to the following MBB  algorithm.
\begin{steps}
  \item Let $b=b(n),1\leq b<n$, be an integer. Denote by $B_t = \{X_t,X_{t+1},\ldots,X_{t+b-1}\}$ the block of length $b$ starting from observation $X_t,\,t=1,2,\ldots,N,$ where $N=n-b+1$ is the number of such blocks available.
  \item Define  i.i.d. integer-valued random variables $I_1, I_2,\ldots, I_k$ having a discrete uniform distribution assigning the probability $1/N$ to each element of the set $\{1,2,\ldots,N\}.$
  \item Let $B_i^*=B_{I_i},\,i=1,2,\ldots,k$, and denote by
$\{X^*_{(i-1)b+1},X^*_{(i-1)b+2},\ldots,X^*_{ib}\}$ the elements of $B_i^*.$ Join the $k$ blocks in the order $B_1^*,B_2^*,\ldots,B_k^*$ together to obtain a new set of functional pseudo observations of length $n$ denoted by $X_1^*,X_2^*,\ldots,X_n^*.$
\end{steps}
\medskip
The above bootstrap algorithm can be potentially applied to approximate the distribution of some statistic $T_n=T(X_1,X_2,\ldots,X_n)$ of interest. For instance, let $T_n=\overline{X}_n$ be the sample mean function of the observed stretch $X_1,X_2,\ldots,X_n,$ i.e., $\overline{X}_n=n^{-1}\sum_{t=1}^{n}X_t.$ We are interested in estimating the distribution of $\sqrt{n}\overline{X}_n$. For this, the bootstrap random variable $\sqrt{n}(\overline{X}^*_n-\E^*(\overline{X}_n^*))$ is used, where $\overline{X}^*_n$ is the mean function of the functional pseudo observations $X_1^*,X_2^*,\ldots,X_n^*,$ i.e., $\overline{X}^*_n=n^{-1}\sum_{t=1}^{n}X^*_t$ and $\E^*(\overline{X}_n^*)$ is the (conditional on the observations $\boldsymbol{\mathrm{X}_n}$) expected value of $\overline{X}^*_n$. Straightforward calculations yield  $$\E^*(\overline{X}_n^*)=\frac{1}{N}\left[\sum_{t=1}^{n}X_t-\sum_{t=1}^{b-1}(1-t/b)(X_t+X_{n-t+1})\right].$$

\medskip
It is known that, under a variety of dependence assumptions on the underlying mean zero stochastic process $\mathds{X},$ it holds true that $\sqrt{n}\overline{X}_n\overset{d}\to \Gamma$ as $n\to\infty,$ where $\Gamma$ denotes a Gaussian process with mean zero and long run covariance operator $2\pi\mathcal{F}_0.$ Furthermore, $\|n\E(\overline{X}_n\otimes\overline{X}_n)-2\pi\mathcal{F}_0\|_{HS}\to 0$ as $n\to\infty$.
Here, $\mathcal{F}_\omega=(2\pi)^{-1}\sum_{h\in\mathbb{Z}}C_he^{-ih\omega},$ $\omega\in\mathbb{R}$, is the so-called spectral density operator of $\mathds{X}$ and $C_h$ denotes the lag $h$ autocovariance operator of $\mathds{X}$, defined by $C_h(\cdot)=\E\langle X_t,\cdot \rangle X_{t+h}$ for any $h\in\mathbb{Z}$; see Panaretos and Tavakoli (2013a,b). 


\medskip

The following theorem establishes validity of the MBB procedure for approximating the distribution of $\sqrt{n}\overline{X}_n$ and for providing a consistent estimator of the long run covariance operator $2\pi\mathcal{F}_0$.
\begin{thm}
\label{thm:CLT}
Suppose that the mean zero stochastic process $\mathds{X}=(X_t, t \in \mathbb{Z})$ satisfies Assumption $1$ and let $X^*_1,X^*_2,\ldots,X^*_n$ be a stretch of pseudo observations generated by the MBB procedure. Assume that the block size $b=b(n)$ satisfies $b^{-1}+bn^{-1/2}=o(1)$ as $n\to\infty.$ Then,  as $n \to \infty,$ 
\begin{enumerate}[label=(\roman*)]
  \item \hfil $\displaystyle d({\cal L}(\sqrt{n}\,(\overline{X}^*_n-\E^*(\overline{X}^*_n)) \mid \boldsymbol{\mathrm{X}_n}), \,{\cal L}(\sqrt{n}\,\overline{X}_n))\to 0, 
\quad\text{in probability,}$\hfill
\end{enumerate}
where $d$ is any metric metrizing weak convergence on $L^2$ and $\mathcal{L}(Z)$ denotes the law of the random element $Z$. Furthermore,
\begin{enumerate}[label=(\roman*)]
 \setcounter{enumi}{1}
  \item \hfil $\displaystyle \|n\E^*(\overline{X}^*_n-\E^*(\overline{X}^*_n))\otimes(\overline{X}^*_n-\E^*(\overline{X}^*_n))-n\E(\overline{X}_n\otimes\overline{X}_n)\|_{HS}=o_P(1),
\quad\text{in probability.}$\hfill
\end{enumerate}
\end{thm}
\medskip

\subsection{\sc The Tapered Block Bootstrap}
\label{ssec:TBB}
The TBB procedure is a modification of the block bootstrap procedure considered in Section~\ref{ssec:MBB} which is obtained by introducing a tapering of the random elements $X_t$. The tapering function down-weights the endpoints of each block $B_i,$ towards zero, i.e., towards the mean function  of $X_t.$ The pseudo observations are then obtained by choosing, with replacement, $k$ appropriately scaled and tapered blocks of length $b$ of centered observations and joining them together.
\par
More precisely, the TBB procedure applied to the functional time series $\boldsymbol{\mathrm{Y}_n}=\{Y_t,\,t=1,2,\ldots,n\}$ stemming from a strictly stationary, $L^2$-values, stochastic process $\mathds{Y}=(Y_t,\,t\in \mathbb{Z})$ with $\E Y_t=0,$ can be described as follows. Let $X_1,X_2,\ldots,X_n$ be the centered observations, i.e., $X_t=Y_t-\overline{Y}_n,\,t=1,2,\ldots,n$, where $\overline{Y}_n=n^{-1}\sum_{t=1}^{n}Y_t$. Furthermore, let $b=b(n)$, $1 \leq b <n$,  be an integer and let $w_n(\cdot),\,n=1,2,\ldots$, be a sequence of so-called data-tapering windows which satisfy the following assumption:
\begin{ass}\label{as:giawtbb}
$w_n(\tau)\in[0,1]$ and $w_n(\tau)=0$ for $\tau\notin\{1,2,\ldots,n\}.$ Furthermore,
\begin{equation}
\label{eq:formwn}
w_n(\tau)=w\left(\dfrac{\tau-0.5}{n}\right),
\end{equation}
where the function $w:\mathbb{R}\to[0,1]$ fulfills the conditions: $(i)$ $w(\tau)\in[0,1]$ for all $\tau\in \mathbb{R}$ with $w(\tau)=0$ if $\tau\notin[0,1]$; $(ii)$ $w(\tau)>0$ for all $\tau$ in a neighbourhood of $1/2$; $(iii)$ $w(\tau)$ is symmetric around $\tau=0.5$; and $(iv)$ $w(\tau)$ is nondecreasing for all $\tau\in[0,1/2].$
\end{ass}
\medskip
Let
$$ \widetilde{B}_i=\left\{w_b(1)\dfrac{b^{1/2}}{\|w_b\|_2}X_i,w_b(2)\dfrac{b^{1/2}} {\|w_b\|_2}X_{i+1},\ldots,w_b(b)\dfrac{b^{1/2}}{\|w_b\|_2}X_{i+b-1} \right\},
$$
be a block of length $b$ starting from $X_t,\,t=1,2,\ldots,N,$ where each centered observation is multiplied by $w_b(\cdot)$ and scaled by $b^{1/2}/\|w_b\|_2,$ $\|w_b\|^2_2=\sum_{i=1}^{b}w_b(i).$ Let $I_1,I_2,\ldots,I_k$  be i.i.d. integers selected from a discrete uniform distribution which assigns probability $1/N$ to each element of the set $\{1,2,\ldots,N\}$. Let $B_i^*=\widetilde{B}_{I_i},\,i=1,2,\ldots,k$, and denote the $i$-th block selected by $\{X^*_{(i-1)b+1},X^*_{(i-1)b+2},\ldots,X^*_{ib}\}$. Join these blocks together in the order $B_1^*,B_2^*,\ldots,B_k^*$ to form the set of TBB pseudo observations $X_1^*,X_2^*,\ldots,X_{n}^*.$\par
Notice that the ``inflation'' factor $b^{1/2}/\|w_b\|_2$ is necessary to compensate for the decrease of the variance of the $X^*_i$'s effected by the shrinking caused by the window $w_b$; see, also, Paparoditis and Politis (2001). Furthermore, the TBB procedure uses the centered time series $X_1,X_2,\ldots,X_{n}$ instead of the original time series $Y_1,Y_2,\ldots,Y_n,$ in order to shrink the end points of the blocks towards zero.\par
To estimate the distribution of $\sqrt{n}\overline{Y}_n$ by means of the TBB procedure, the bootstrap random variable $\sqrt{n}(\overline{X}^*_n-\E^*(\overline{X}_n^*))$ is used, where $\overline{X}^*_n=n^{-1}\sum_{t=1}^{n}X^*_t$ and $\E^*(\overline{X}_n^*)$ is the (conditional on the observations $\boldsymbol{\mathrm{Y}_n}$) expected value of $\overline{X}^*_n.$ Straightforward calculations yield $$\E^*(\overline{X}_n^*)=\dfrac{1}{N}\dfrac{\|w_b\|_1}{\|w_b\|_2}\left[\sum_{t=1}^{n} X_t-\sum_{t=1}^{b-1}\left(1-\dfrac{\sum_{s=1}^{t}w_b(s)}{\|w_b\|_1}\right)X_t-\sum_{j=1}^{b-1}\left(1-\dfrac{\sum_{t=b-j+1}^{b}w_b(t)}{\|w_b\|_1}\right)X_{n-j+1}\right].$$

\medskip
The following theorem establishes validity of the TBB procedure for approximating the distribution of $\sqrt{n}\overline{Y}_n$ and for providing a consistent estimator of the long run covariance operator $2\pi\mathcal{F}_0$.

\begin{thm}
\label{thm:CLTtbb}
Suppose that the mean zero stochastic  process $\mathds{Y}$ satisfies Assumption~\ref{as:L2} and let $w_n(\cdot),\,n=1,2,\ldots,$ be a sequence of data-tapering windows satisfying Assumption~\ref{as:giawtbb}. Furthermore, let $X^*_t, t=1,2,\ldots,n,$ be a stretch of pseudo observations generated by the TBB procedure. Assume that the block size $b=b(n)$ satisfies
$b^{-1}+bn^{-1/2}=o(1)$ as $n\to\infty$. Then,  as $n\to\infty$, 
\begin{enumerate}[label=(\roman*)]
 \item \hfil $\displaystyle d({\cal L}(\sqrt{n}\,(\overline{X}^*_n-\E^*(\overline{X}^*_n)) \mid \boldsymbol{\mathrm{Y}_n}),\, {\cal L}(\sqrt{n}\,\overline{Y}_n))\to 0,\quad\text{in probability,}$\hfill
\end{enumerate}
where $d$ is any metric metrizing weak convergence on $L^2$,
and
\begin{enumerate}[label=(\roman*)]
 \setcounter{enumi}{1}
  \item \hfil $\displaystyle \|n\E^*(\overline{X}^*_n-\E^*(\overline{X}^*_n))\otimes(\overline{X}^*_n-\E^*(\overline{X}^*_n))-n\E(\overline{Y}_n\otimes\overline{Y}_n)\|_{HS}=o_P(1),
\quad\text{in probability.}$\hfill
\end{enumerate}
\end{thm}

\begin{remark}
The asymptotic validity of the MBB and TBB procedures established in Theorem~\ref{thm:CLT}  and Theorem~\ref{thm:CLTtbb}, respectively, can be extended to cover  also the case where  maps  $\phi : L ^2\rightarrow D$ of  the sample means $\overline{X}_n$ (in the MBB case) and $ \overline{Y}_n$ (in the TBB case) are  considered, when $D$ is a normed space. For instance,  such a result follows  as  an application of  a version of the delta-method appropriate for the bootstrap and  for maps $\phi$ which are Hadamard differentiable at $0$ tangentially to a subspace $D_0$ of $D$ (see Theorem 3.9.11 of van der Vaart and Wellner (1996)).  Extensions of such results to almost surely convergence and for  other types of differentiable maps, like for instance Fr\'echet differentiable functionals (see Theorem  3.9.15 of van de Vaart and Wellner (1996))  or quasi-Hadamard differentiable functionals (see Theorem 3.1 of Beutner and Z\"ahe (2016)),  are more involved since they depend on the particular map $\phi$ and  the verification of  some technical  conditions. 
\end{remark}

\section{\sc Bootstrap-Based Testing of the Equality of Mean Functions}
\label{sec:testing}
Among different applications, the MBB and TBB procedures can be also used to perform a test of the equality of mean functions between several independent samples of a functional time series. In this case, both block bootstrap procedures have to be implemented in such a way that the pseudo observations $X_1^*,X_2^*,\ldots,X_n^*$ generated, satisfy the null hypothesis of interest. 
\subsection{\sc The set-up}
\label{ssec:algorithm}
Consider $K$ independent functional time series $\boldsymbol{\mathrm{X_M}}=\{X_{i,t},\,i=1,2\ldots,K,\,t=1,2,\ldots,n_i\}$, each one of which satisfies
\begin{equation}
\label{model:testing}
X_{i,t}=\mu_i+\varepsilon_{i,t},\,\,t=1,2,\ldots,n_i,
\end{equation}
where, for each $i\in\{1,2,\ldots,K\},\,\{\varepsilon_{i,t},t\in\mathbb{Z}\}$ is a $L^2$-$m$-approximable functional process and $n_i$ denotes the length of the $i$-th time series. Let $M = \sum_{i=1}^Kn_i$ be the total number of observations and note that $\mu_i(\tau),\,\tau\in \mathcal{I},$ is the mean function of the $i$-th functional time series, $i=1,2,\ldots,K.$ The null hypothesis of interest is then, $$H_0:\mu_1=\mu_2=\ldots=\mu_K $$ and the alternative hypothesis
$$ H_1:\exists\, k,m \in\{1,2,\ldots,K\}\; \text{with}\; k\neq m\; \text{such that}\; \mu_k\neq \mu_m .$$
Notice that the above equality is in $L^2$, i.e., $\mu_k=\mu_m$ means that $\|\mu_m-\mu_k\|=0$ whereas $\mu_k\neq\mu_m$ that $\|\mu_m-\mu_k\|>0.$

\subsection{\sc Block Bootstrap-based testing}
The aim is to generate a set of functional pseudo observations $\boldsymbol{\mathrm{X^*_M}}=X^*_{i,t},\,i=1,2\ldots,K$, $t=1,2,\ldots,n_i,$ using either the MBB procedure or the TBB procedure in such a way that $H_0$ is satisfied. These bootstrap pseudo-time series can then be used to estimate the distribution of some test statistic $T_M$ of interest which is applied to test $H_0.$ Toward this, the distribution of $T^*_M$ is used as an estimator of the distribution of $T_M$, where $T^*_M$ is the same statistical functional as $T_M$ but calculated using the bootstrap functional pseudo-time series $\boldsymbol{\mathrm{X^*_M}}.$ \par
To implement the MBB procedure for testing the null hypothesis of interest, assume, without loss of generality, that the test statistic $T_M$ rejects the null hypothesis when $T_M > d_{M,\alpha},$ where, for $\alpha\in(0,1)$,  $d_{M,\alpha}$ denotes the upper $\alpha$-percentage point of the distribution of $T_M$ under $H_0.$
The MBB-based testing procedure goes then as follows:
\begin{steps}
\item Calculate the sample mean functions in each population and the pooled mean function, i.e., calculate
$\overline{X}_{i,n_i}=(1/n_i)\sum_{t=1}^{n_i}X_{i,t}$, for $i=1,2\ldots,K$, and $\overline{X}_M=(1/M)\sum_{i=1}^{K}\sum_{t=1}^{n_i}X_{i,t}$,  and obtain the residual functions
 in each population, i.e., calculate $\hat{\varepsilon}_{i,t}=X_{i,t}-\overline{X}_{i,n_i}$, for $t=1,2,\ldots,n_i$;  $i=1,2\ldots,K$.
\item For $i=1,2,\ldots,K$, let $b_i=b_i(n)\in\{1,2,\ldots,n-1\}$ be the block size for the $i$-th functional time series and divide $\{\hat{\varepsilon}_{i,t},\,t=1,2,\ldots,n_i\}$ into $N_i=n_i-b_i+1$ overlapping blocks of length $b_i$, say, $B_{i,1},B_{i,2},\ldots,B_{i,N_i}.$ Calculate the sample mean of the $\xi$-th observations of the blocks $B_{i,1},B_{i,2},\ldots,B_{i,N_i}$, i.e.,
$\overline{\varepsilon}_{i,\xi}=(1/N_i)\sum_{t=1}^{N_i}\hat{\varepsilon}_{i,\xi+t-1}$, for $\xi=1,2,\ldots, b_i$.
\item For simplicity assume that $n_i=k_ib_i$ and for $i=1,2,\ldots,K$, let $q^i_1,q^i_2,\ldots, q^i_{k_i}$ be i.i.d. integers selected from a discrete probability distribution which assigns the probability $1/N_i$ to each element of the set $\{1,2,\ldots,N_i\}.$ Generate bootstrap functional pseudo observations $X^*_{i,t},\,t=1,2,\ldots,n_i,\,i=1,2,\ldots,K$, as
$
X^*_{i,t}=\overline{X}_M+\varepsilon^*_{i,t},
$
where
\begin{equation}
\label{eq:epsilostar}
\varepsilon^*_{i,\xi+(s-1)b_i}=\hat{\varepsilon}_{i,q^i_s+\xi-1}-\overline{\varepsilon}_{i,\xi},\,s=1,2,\ldots,k_i,\,\xi=1,2,\ldots,b_i.
\end{equation}
\item Let $T_M^*$ be the same statistic as $T_M$ but calculated using the bootstrap functional pseudo-time series $X^*_{i,t},\,t=1,2,\ldots,n_i$, $i=1,2,\ldots,K$. Denote by $D^*_{M,T}$
the distribution of $T^*_M$ given $\boldsymbol{\mathrm{X_M}}.$ For $\alpha \in (0,1),$ reject the null hypothesis $H_0$ if
$T_M > d_{M,\alpha}^*$,
where $d_{M,\alpha}^*$ denotes the upper $\alpha$-percentage point of the distribution of $T_M^*$, i.e., $\mathbb{P}(T_M^*>d_{M,\alpha}^*)=\alpha$.
\end{steps}
Note that the distribution $D^*_{M,T}$ can be evaluated by Monte-Carlo.\par
\medskip
To motivate the centering used in Step~$3,$ denote, for $i=1,2,\ldots,K$, by $e^*_{i,t},\,t=1,2,\ldots,n_i,$ the pseudo observations generated by applying the MBB procedure, described in Section~\ref{ssec:MBB}, directly to the residuals time series $\hat{\varepsilon}_{i,t},\,t=1,2,\ldots,n_i.$ Note that the $e^*_{i,t}$'s differ from the $\varepsilon^*_{i,t}$'s used in~\eqref{eq:epsilostar} by the fact that the later are obtained after centering. The sample mean $\overline{\varepsilon}_{i,\xi}$, $i=1,2,\ldots,K,$ $\xi=1,2,\ldots,b_i$, calculated in Step~2, is the (conditional on $\boldsymbol{\mathrm{X_M}})$ expected value of the pseudo observations $e^*_{i,t},\,t=1,2,\ldots,n_i,$ where $t=\xi\pmod{b_i}.$ Furthermore, for $i=1,2,\ldots,K,$ we generate the $\varepsilon_{i,t}^*$'s, $t=1,2,\ldots,n_i$, by subtracting $\overline{\varepsilon}_{i,\xi}$ from $e^*_{i,sb+\xi},\,\xi=1,2,\ldots,b,\,s=1,2,\ldots k_i.$
This is done in order for the (conditional on $\boldsymbol{\mathrm{X_M}}$) expected value of $\varepsilon^*_{i,t}$ to be zero. In this way, the generated set of pseudo time series  $X^*_{i,t},\,t=1,2,\ldots,n_i,\,i=1,2,\ldots,K,$ satisfy the null hypothesis $H_0$. In particular, given $\boldsymbol{\mathrm{X_M}}=\{X_{i,t},\,i=1,2\ldots,K,\,t=1,2,\ldots,n_i\}$, we have
$$
\mathbb{E}^*(X^*_{i,\xi+(s-1)b_i}) 
=\overline{X}_M+\dfrac{1}{N_i}\sum_{t=1}^{N_i}[\hat{\varepsilon}_{i,t+\xi-1}-\overline{\varepsilon}_{i,\xi}]
=\overline{X}_M,\label{eq:mesialgorithmouMBB}
$$
for $i=1,2\ldots,K,$ $\xi=1,2,\ldots,b$ and $s=1,2,\ldots,k_i.$
That is, conditional on $\boldsymbol{\mathrm{X_M}},$ the mean function of each functional pseudo-time series $X_{i,1}^*,X_{i,2}^*,\ldots,X_{i,n}^*,$ $i=1,2\ldots,K,$ is identical in each population and equal to the pooled sample mean function $\overline{X}_M$.\par 
An algorithm based on the TBB procedure for testing the same pair of hypotheses can also be implemented by modifying appropriate the MBB-based testing algorithm. In particular, we replace Step~$2$ and Step~$3$ of this algorithm by the following steps:
\begin{steps}\addtocounter{step}{1}
\item For $i=1,2,\ldots,K$, let $b_i=b_i(n)\in\{1,2,\ldots,n-1\}$ be the block size for the $i$-th functional time series and $N_i=n_i-b_i+1.$ Let also $\{\hat{\epsilon}_{i,t},\,t=1,2,\ldots,n_i\}$ be the centered values of $\{\hat{\varepsilon}_{i,t},\,t=1,2,\ldots,n_i\},$ i.e., $\hat{\epsilon}_{i,t}=\hat{\varepsilon}_{i,t}-\overline{\varepsilon}_{i}$, where $\overline{\varepsilon}_{i}=(1/n_i)\sum_{t=1}^{n_i}\hat{\varepsilon}_{i,t}.$ Also, let $w_{n_i}(\cdot),\,n_i=1,2,\ldots,$ be a sequence of data-tapering windows satisfying Assumption~\ref{as:giawtbb}. Now, for $t=1,2,\ldots,N_i$, let
$$
\widetilde{B}_{i,t}=\left\{w_{b_i}(1)\dfrac{b_i^{1/2}}{\|w_{b_i}\|_2}\hat\epsilon_{i,t},w_{b_i}(2)\dfrac{b_i^{1/2}} {\|w_{b_i}\|_2}\hat\epsilon_{i,t+1},\ldots,w_{b_i}(b_i)\dfrac{b_i^{1/2}}{\|w_{b_i}\|_2}\hat\epsilon_{i,t+b_i-1} \right\},\quad i=1,2,\ldots,K,
$$
where $\|w_{b_i}\|^2_2=\sum_{i=1}^{b_i}w_{b_i}(i).$ Here, $\widetilde{B}_{i,t}$ denotes the tapered block of $\hat{\epsilon}_{i,t}$'s of length $b_i$ starting from $\hat\epsilon_{i,t}.$ Furthermore, for $i=1,2,\ldots,K$, calculate the sample mean of the $\xi$th observations of the blocks $\widetilde{B}_{i,1},\widetilde{B}_{i,2},\ldots,\widetilde{B}_{i,N_i}$, i.e.,
$$ \bar{\epsilon}_{i,\xi}=\dfrac{1}{N_i}\sum_{t=1}^{N_i}w_{b_i}(\xi)\dfrac{b_i^{1/2}}{\|w_{b_i}\|_2}\hat{\epsilon}_{i,\xi+t-1},\,\,\xi=1,2,\ldots,b_i. $$
\item  For $i=1,2,\ldots,K,$ let $q^i_1, q^i_2,\ldots, q^i_{k_i}$ be i.i.d. integers selected from a discrete probability distribution which assigns the probability $1/N_i$ to each $t\in\{1,2,\ldots,N_i\}.$ Generate bootstrap functional pseudo-observations $X^+_{i,t},\,i=1,2,\ldots,K,\,t=1,2,\ldots,n_i$ according to
$
X^+_{i,t}=\overline{X}_M+\epsilon^+_{i,t},
$
where
$$
\epsilon^+_{i,\xi+(s-1)b_i}=w_b(\xi)\dfrac{b_i^{1/2}}{\|w_{b_i}\|_2}\hat{\epsilon}_{i,k^i_s+\xi-1}-\bar{\epsilon}_{i,\xi},\,s=1,2,\ldots,k_i,\,\xi=1,2,\ldots,b_i.
$$
\end{steps}
As in the case of the MBB-based testing, the generation of $\epsilon^+_{i,t},\,t=1,2,\ldots,n_i,\,i=1,2,\ldots,K,$ ensures that the functional pseudo-time series $X^+_{i,t},\,t=1,2,\ldots,n_i,\,i=1,2,\ldots,K,$ satisfy $H_0,$ that is, given $\boldsymbol{\mathrm{X_M}}=\{X_{i,t},\,i=1,2\ldots,K,\,t=1,2,\ldots,n_i\}$, we have that
$
\mathbb{E}^+(X^+_{i,t})=\overline{X}_M.
$
\subsection{\sc Bootstrap Validity}
\label{sec:consistency}
Notice that, since the proposed block bootstrap-based methodologies are not designed having any particular test statistic in mind, they can be potentially applied to a wide range of test statistics. To prove validity of the proposed block bootstrap-based testing procedures, however, a particular test statistic has to be considered. 
For instance, one such  test statistic is the fully functional test statistic proposed by Horv\'{a}th {\em et al.} (2013) for the case of $K=2$ populations. Let $X_{i,t},\,i=1,2,\,t=1,2,\ldots,n_i,$ be two independent samples of curves, satisfying model~\eqref{model:testing}. 
For $i\in\{1,2\}$ and for $(u,v)\in[0,1]^2$, denote by $c_i(u,v)$ the kernels
of the long run covariance operators $2\pi\mathcal{F}^{(i)}_0$, 
given by
$
c_i(u,v)=\E[\varepsilon_{i,0}(u)\varepsilon_{i,0}(v)]+\sum_{j\geq1}\E[\varepsilon_{i,0}(u)\varepsilon_{i,j}(v)]+\sum_{j\geq1}\E[\varepsilon_{i,0}(v)\varepsilon_{i,j}(u)].
$
The test statistic considered in  Horv\'{a}th {\em et al.} (2013), evaluates  then the $L^2$-distance of the two sample mean functions $\overline{X}_{i,n_i}=n_i^{-1}\sum_{t=1}^{n_i}X_{i,t},\,i=1,2,$ and it is given by
$$
U_M=\dfrac{n_1n_2}{M}\int(\overline{X}_{1,n_1}(\tau)-\overline{X}_{2,n_2}(\tau))^2\,\mathrm{d}\tau,
$$
where $M=n_1+n_2.$ Horv\'{a}th {\em et al.} (2013) proved that if $\min\{n_1,n_2\}\to\infty$ and $n_1/M\to\theta\in(0,1)$ then, under $H_0,$ $U_M$ converges weakly to $\int\Gamma^2(\tau)\,\mathrm{d}\tau$, where $\{\Gamma(\tau):\,\tau\in\mathcal{I}\}$ is a Gaussian process satisfying $\mathbb{E}(\Gamma(\tau))=0$ and $\mathbb{E}(\Gamma(u)\Gamma(v))=(1-\theta)c_1(u,v)+\theta c_2(u,v).$
Notice that calculation of critical values of the above test requires estimation of the distribution of $\int\Gamma^2(\tau)\,\mathrm{d}\tau$ which is a difficult task.

Although the test statist $U_M$ is quite appealing because it is fully functional, its limiting distribution is difficult to implement which demonstrates  the importance of the bootstrap. To investigate the consistency properties of the bootstrap, we first establish a general result which allows for the consideration of different test statistics that can be expressed as functionals  of the basic deviation process  
\begin{equation} \label{eq.devproc}
  \Big\{ \sqrt{\dfrac{n_1n_2}{M}}\big(\overline{X}_{1,n_1}(\tau)-\overline{X}_{2,n_2}(\tau) \big), \tau\in {\mathcal I}\Big\}.
\end{equation} 

\begin{thm}
\label{thm:basic-consistent}
Let Assumption~\ref{as:L2} be satisfied. Assume that $\min\{n_1,n_2\}\to\infty,$ $n_1/M\to\theta\in(0,1)$ and that, for $i\in\{1,2\}$, the block size $b_i=b_i(n)$ fulfills  $b_i^{-1}+b_in_i^{-1/2}=o(1),$ as $n_i\to\infty$. Then, conditional on $\boldsymbol{\mathrm{X_M}}$, as $n_i\to\infty$, 
\begin{enumerate}[label=(\roman*)]
\item  $ \sqrt{\dfrac{n_1n_2}{M}}\big(\overline{X}^*_{1,n_1}-\overline{X}^*_{2,n_2} \big) 
\Rightarrow  \Gamma$, in probability,
\end{enumerate}
and, if additionally Assumption~\ref{as:giawtbb} is satisfied, 
\begin{enumerate}[label=(\roman*)]
 \setcounter{enumi}{1}
\item 
$ \sqrt{\dfrac{n_1n_2}{M}}\big(\overline{X}^+_{1,n_1}-\overline{X}^+_{2,n_2} \big) \Rightarrow  \Gamma$, in probability.\\
\end{enumerate}
Here, $ \Rightarrow  $ denotes weak convergence in $ L^2$. 
\end{thm}

By Theorem~\ref{thm:basic-consistent} and  the continuous mapping theorem,  the suggested block bootstrap-based testing procedures can be successfully applied to consistently estimate the distribution of 
any test statistic of interest which is a continuous function of the basic deviation process (\ref{eq.devproc}).  We elaborate on some examples. Below, $P_{H_0}(Z\leq \cdot)$ denotes the distribution function of the random variable $Z$ when $H_0$ is true.
\medskip

Consider for  instance  the test statistic $U_M$.  
Let
$$
U^*_M=\dfrac{n_1n_2}{M}\int(\overline{X}^*_{1,n_1}(\tau)-\overline{X}^*_{2,n_2}(\tau))^2\,\mathrm{d}\tau
$$
and
$$
U^+_M=\dfrac{n_1n_2}{M}\int(\overline{X}^+_{1,n_1}(\tau)-\overline{X}^+_{2,n_2}(\tau))^2\,\mathrm{d}\tau,
$$
where
$
\overline{X}_{i,n_i}^*=(1/n_i)\sum_{t=1}^{n_i}X_{i,t}^*$ and $
\overline{X}_{i,n_i}^+=\dfrac{1}{n_i}\sum_{t=1}^{n_i}X_{i,t}^+$, $i=1,2$.
We then have the following result.
\begin{corollary}
\label{thm:consistent}
Let the assumptions of Theorem~\ref{thm:consistent} be satisfied. Then, 
\begin{enumerate}[label=(\roman*)]
\item $
\sup_{x\in \mathbb{R}} \bigl| P(U_M^*\leq x \mid \boldsymbol{\mathrm{X_M}})-P_{H_0}(U_M\leq x) \bigr|\to 0,\,\,\,\,\text{in probability}, 
$
and
\item $
\sup_{x\in \mathbb{R}} \bigl| P(U_M^+\leq x \mid \boldsymbol{\mathrm{X_M}})-P_{H_0}(U_M\leq x) \bigr|\to 0,\,\,\,\,\text{in probability}.
$
\end{enumerate}
\end{corollary}


\begin{remark}
If the following type of  one-sided tests  $H_0: \mu_1=\mu_2$ versus $H_1: \mu_1 >\mu_2$ or $H'_1: \mu_1 <\mu_2$ 
is of interest (where $\mu_1=\mu_2$ (resp $\mu_1 >\mu_2$ or $\mu_1 <\mu_2$)  means $\mu_1(\tau)=\mu_2(\tau)$ (resp $\mu_1(\tau) >\mu_2(\tau)$ or $\mu_1(\tau) <\mu_2(\tau)$) for all $\tau \in {\cal I}$), then the following test statistic
$$
\widetilde{U}_M=\sqrt{\dfrac{n_1n_2}{M}}\int(\overline{X}_{1,n_1}(\tau)-\overline{X}_{2,n_2}(\tau))\,\mathrm{d}\tau
$$
can be used. In this case, $H_0$ is rejected if $\widetilde{U}_M > \widetilde{d}_{M,\alpha}$ or $\widetilde{U}_M < \widetilde{d}_{M,\alpha}$, respectively, with $\widetilde{d}_{M,\alpha}$ the upper $\alpha$-percentage point of the distribution of $\widetilde{U}_M$ under $H_0.$ Consistent estimators of this distribution can be also obtained using the block bootstrap procedures discussed. In particular, the following results can be established:
\begin{enumerate}[label=(\roman*)]
\item $
\sup_{x\in \mathbb{R}} \bigl| P(\tilde{U}_M^*\leq x \mid \boldsymbol{\mathrm{X_M}})-P_{H_0}(\tilde{U}_M\leq x) \bigr|\to 0,\,\,\,\,\text{in probability},\,  and
$
\item $
\sup_{x\in \mathbb{R}} \bigl| P(\tilde{U}_M^+\leq x \mid \boldsymbol{\mathrm{X_M}})-P_{H_0}(\tilde{U}_M\leq x) \bigr|\to 0,\,\,\,\,\text{in probability},
$
\end{enumerate}
where
$$
\widetilde{U}^*_M=\sqrt{\dfrac{n_1n_2}{M}}\int(\overline{X}^*_{1,n_1}(\tau)-\overline{X}^*_{2,n_2}(\tau))\,\mathrm{d}\tau
$$
and
$$
\widetilde{U}^+_M=\sqrt{\dfrac{n_1n_2}{M}}\int(\overline{X}^+_{1,n_1}(\tau)-\overline{X}^+_{2,n_2}(\tau))\,\mathrm{d}\tau.
$$
To elaborate, notice that using Theorem~$1$ of Horv\'{a}th {\em et al.} (2013), we get, as $n_1,n_2\to\infty,$ that
$$
\left(\dfrac{1}{\sqrt{n_1}}\sum_{j=1}^{n_1}(X_{1,j}-\mu_1),\dfrac{1}{\sqrt{n_2}}\sum_{j=1}^{n_2}(X_{2,j}-\mu_2)\right)
\Rightarrow (\Gamma_1,\Gamma_2),
$$
where $\Gamma_1$ and $\Gamma_2$ are two independent Gaussian random elements in $L^2$ with mean zero and covariance operators $C_1$ and $C_2$ with kernels $c_1(\cdot,\cdot)$ and $c_2(\cdot,\cdot),$ respectively. Under $H_0,$ and for $\tilde{\mu}=\mu_1=\mu_2$ the common mean of the two populations, we have
$$
\sqrt{\dfrac{n_1n_2}{M}}(\overline{X}_{1,n_1}(\tau)-\overline{X}_{2,n_2}(\tau))=
\sqrt{\dfrac{n_2}{M}}\dfrac{1}{\sqrt{n_1}}\sum_{t=1}^{n_1}(X_{1,t}-\tilde{\mu})-\sqrt{\dfrac{n_1}{M}}\dfrac{1}{\sqrt{n_2}}\sum_{t=1}^{n_2}(X_{2,t}-\tilde{\mu}),
$$
which implies, for $n_1,n_2\to\infty$ and $n_1/M\to\theta,$ that
$\widetilde{U}_M\overset{d}\to\int\Gamma(\tau)\,\mathrm{d}\tau,$ where $\Gamma(\tau)=\sqrt{1-\theta}\Gamma_1(\tau)-\sqrt{\theta}\Gamma_2(\tau)$, $\tau \in {\cal I}$. Now, working along the same lines as in the proof of Theorem~\ref{thm:consistent}, it can be easily shown that $\widetilde{U}^*_M$ and $\widetilde{U}^+_M$ converges weakly to the same limit $\int\Gamma(\tau)\,\mathrm{d}\tau.$
\end{remark}

Another interesting class of test statistics for which Theorem~\ref{thm:basic-consistent} allows for a  successful application of 
the  suggested block bootstrap-based testing  procedures, is the class of so-called projection-based tests. To elaborate, let $\{\varphi_1, \varphi_2, \ldots, \varphi_p\}$ be a set of $p$ orthonormal functions in $L_2$. A common choice is  to let $\varphi_j$ be  the  orthonormalized eigenfunctions  corresponding to the $p$  largest eigenvalues of the covariance 
operator of the stochastic process $\{\Gamma(\tau)=\sqrt{1-\theta}\Gamma_1(\tau)-\sqrt{\theta}\Gamma_2(\tau),\, \tau \in \cal{I}\}$, which are assumed to be distinct and strictly positive. A  test statistic $S_{p,M}$ can then be considered which is defined as 
\[  S_{p,M}= \frac{n_1 n_2}{M}\sum_{k=1}^p\langle \overline{X}_{1,n_1}-\overline{X}_{2,n_2}, \widehat{\varphi}_k \rangle^2,\]
and where $ \widehat{\varphi}_k $ are  estimators of $ \varphi_k$; see for instance Horv\'{a}th {\em et al.} (2013) where studentized versions of $\langle \overline{X}_{1,n_1}-\overline{X}_{2,n_2}, \widehat{\varphi}_k \rangle$ have also been used.
\medskip

The following result establishes consistency of the suggested block bootstrap methods also for this class of test statistics. 

\begin{corollary}
\label{thm:consistent2}
Let the assumptions of Theorem~\ref{thm:consistent} be satisfied and assume that the $p$ largest eigenvalues of  the covariance operator of the stochastic process $\{\Gamma(\tau)=\sqrt{1-\theta}\Gamma_1(\tau)-\sqrt{\theta}\Gamma_2(\tau),\, \tau \in \cal{I}\}$ are distinct and positive. Let  $\varphi_k$, $ k=1,2, \ldots, p$, be the orthonormalized eigenfunctions corresponding to these eigenvalues and 
let $ \widetilde{\varphi}_k$ and $ \widehat{\varphi}_k$ be estimators of $ \varphi_k$ satisfying $ \max_{1\leq k\leq p}\|\widetilde{\varphi}_k-\widetilde{c}_k\varphi_k\| \stackrel{P}{\rightarrow} 0$ and  
$\max_{1\leq k\leq p}\|\widehat{\varphi}_k-\widehat{c}_k\varphi_k\| \stackrel{P}{\rightarrow} 0$, where $ \widetilde{c}_k= sign\big(\langle \widetilde{\varphi}_k,\varphi_k\rangle)$ and 
$ \widehat{c}_k= sign\big(\langle \widehat{\varphi}_k,\varphi_k\rangle) $. Then,  
\begin{enumerate}[label=(\roman*)]
\item $
\sup_{x\in \mathbb{R}} \bigl| P(S_{p,M}^*\leq x \mid \boldsymbol{\mathrm{X_M}})-P_{H_0}(S_{p,M}\leq x) \bigr|\to 0,\,\,\,\,\text{in probability}, \, and
$
\item $
\sup_{x\in \mathbb{R}} \bigl| P(S_{p,M}^+\leq x \mid \boldsymbol{\mathrm{X_M}})-P_{H_0}(S_{p,M}\leq x) \bigr|\to 0,\,\,\,\,\text{in probability},
$
\end{enumerate}
where $S_{p,M}^*= (n_1 n_2/ M)\sum_{k=1}^p\langle \overline{X}^*_{1,n_1}-\overline{X}^*_{2,n_2}, \widetilde{\varphi}_k \rangle^2$
 and $ S_{p,M}^+=(n_1 n_2/ M)\sum_{k=1}^p\langle \overline{X}^+_{1,n_1}-\overline{X}^+_{2,n_2}, \widetilde{\varphi}_k \rangle^2$.
\end{corollary} 

\begin{remark}
In Corollary~\ref{thm:consistent2}, we allow for $ \widetilde{\varphi}_k$ to be a different estimator of $ \varphi_k$ than $\widehat{\varphi}_k$, where the later is  used in the test statistic $ S_{p,M}$. For instance, $ \widetilde{\varphi}_k$ could be the same  estimator as $ \widehat{\varphi}_k$ but based on the the bootstrap pseudo observations $ X^\ast_{i,t}$, $i=1,2, \ldots, k$ and $t=1,2, \ldots, n_i$, respectively,  $ X^+_{i,t}$, $i=1,2, \ldots, k$ and $t=1,2, \ldots, n_i$. This will allow for the bootstrap statistics $ S^\ast_{p,M}$, respectively  $S_{p,M}^+$, 
to also imitate the effect of the estimation  error of the unknown eigenfunctions $\varphi_k$ on the distribution of $S_{p,M}$. Clearly, a simple and computationally easier  
alternative will be to set   
 $ \widetilde{\varphi}_k=\widehat{\varphi}_k$. 
\end{remark}

\begin{remark}
If the alternative hypothesis $H_1$ is true, then under the same assumptions as in Theorem~$4$ of Horv\'{a}th {\em et al.} (2013), we get that $U_M\overset{P}\to\infty.$ Furthermore, under the same assumptions as in Theorem~$6$ of Horv\'{a}th {\em et al.} (2013), we get that $S_{p,M} \overset{P}\to\infty$  provided that $\langle \mu_1-\mu_2,  \varphi_k \rangle \neq 0$ for at least one $1 \leq k \leq p$.
Corollary~\ref{thm:consistent} and Corollary \ref{thm:consistent2} (together with Slutsky's theorem) imply then, respectively, the consistency of the test $U_M$ using the bootstrap  critical values obtained from the distributions of $U_M^*$ and $U_M^+$, and of the test $S_{p,M}$ using the bootstrap  critical values obtained from the distributions of $S_{p,M}^*$ and $S_{p_M}^+$.
\end{remark}

\section{\sc Numerical Examples}
We generated functional time series stemming from a first order functional autoregressive model (FAR(1))
\begin{equation}
\label{model:FAR1}
\varepsilon_t(u)=\int\psi(u,v)\varepsilon_{t-1}(v)\,\mathrm{d}v+B_t(u), \quad u \in [0,1],
\end{equation}
(see also Horv\'{a}th {\em et al.} (2013)), and from a first order functional moving average model (FMA(1)),
\begin{equation}
\label{model:FMA1}
\varepsilon_t(u)=\int\psi(u,v)B_{t-1}(v)\,\mathrm{d}v+B_t(u), \quad u \in [0,1].
\end{equation}
For both models, the kernel function $\psi(\cdot,\cdot)$ is defined by
\begin{equation}
\label{eqn:psikernelF}
\psi(u,v)=\dfrac{\mathrm{e}^{-(u^2+v^2)/2}}{4\int\mathrm{e}^{-t^2}\mathrm{d}t},\,\,(u,v)\in[0,1]^2,
\end{equation}
and the $B_t$'s are i.i.d. Brownian bridges. All curves were approximated using $T=21$ equidistant points $\tau_1,\tau_2,\ldots,\tau_{21}$ in the unit interval $\mathcal{I}$ and transformed into functional objects using the Fourier basis with $21$ basis functions (see Section 3 of the supplementary material for additional simulations with a larger $T$).\par

Implementation of the MBB and TBB procedures require the selection of the block size $b.$ As it has been shown in Theorem~\ref{thm:CLT} and Theorem~\ref{thm:CLTtbb},
$n\E^*[(\overline{X}_n^*-\E^*(\overline{X}_n^*))\otimes(\overline{X}_n^*-\E^*(\overline{X}_n^*))]$ is a consistent estimator of $2\pi\mathcal{F}_0,$ with  kernel
$$
c_N(u,v)=\dfrac{1}{N}\sum_{i=1}^{n} X_i(u)X_i(v)+
\sum_{h=1}^{b-1}\left(1-\dfrac{h}{b}\right)
\dfrac{1}{N}\sum_{i=1}^{n-h}\Big[X_i(u)X_{i+h}(v)+X_{i+h}(u)X_i(v)\Big]+o_p(1),
$$
in the MBB case, and
$$
\tilde{c}_N(u,v)=\dfrac{1}{N}\sum_{i=1}^{n} Y_i(u)Y_i(v)+
\sum_{h=1}^{b-1}\dfrac{\mathcal{W}_h}{\|w_b\|_2^2}
\dfrac{1}{N}\sum_{i=1}^{n-h}[Y_i(u)Y_{i+h}(v)+Y_{i+h}(u)Y_i(v)]+o_p(1),
$$
in the TBB case, with $\mathcal{W}_h=\sum_{i=1}^{b-h}w_b(i)w_b(i+h),\,h=0,1,\ldots,b-1$. Now, 
$c_N$ and $\tilde{c}_N$ can be considered as lag-window estimators of the kernel $c(u,v)=\sum_{i=-\infty}^{\infty}\E[X_0(u)X_i(v)],$ using the Bartlett window with ``truncation lag" $b$ in the MBB case and using the same ``truncation lag" with the window function $W=\mathcal{W}_h/\|w_b\|$, in the TBB case.
The above observations suggest that the problem of choosing the block size $b$ can be considered as a problem of choosing the truncation lag of a lag window estimator of the function $c(u,v).$ Choice of the truncation lag in the functional context has been recently discussed in Horv\'{a}th {\em et al.} (2016) and Rice and Shang (2016).
Although different procedures to select the ``truncation lag" have been proposed in the aforementioned papers, we found the simple rule of setting $b_i=\left\lceil n_i^{1/3}\right\rceil$, where $\lceil x \rceil$ is the least integer greater than or equal to $x$, quite effective in our numerical examples. In the following, we denote by $b^*$ this choice of $b$, which is used together with some other values of $b_i.$ A simulation study has been first conducted in order to investigate the finite sample performance of the MBB and TBB procedures. For this, the problem of estimating the standard deviation function of the normalized sample mean $\sqrt{n}\overline{X}_n(\tau)$, i.e., of $\sqrt{c(\tau,\tau)}$ for different values of $\tau \in [0,1]$ has been considered. The results obtained using both block bootstrap procedures have also been compared with those using the stationary bootstrap (SB).
Realizations of length $n=100$ and $n=500$ from the functional time series models~\eqref{model:FAR1} and~\eqref{model:FMA1} have been used. The results obtained are presented and discussed in Section 2 of the supplementary material.  Furthermore, Table 1 of the supplementary material presents results comparing the performance of projections-based tests when asymptotic and bootstrap approximations  are used to obtain the critical values of  the tests. 

\subsection{\sc Testing equality of mean functions}
\label{ssec:simulations1b}
We investigate the size and power performance of the tests considered in Section~\ref{sec:consistency}. As can be seen in Section 2 of the supplementary material, the TBB estimators perform best in our simulations. For this reason, we concentrate in this section, on tests based on TBB critical values only. Two sample sizes  $n_1=n_2=100$ and $n_1=n_2=200$ as well as a range of block sizes $b=b_1=b_2,$ are considered. The tests have been applied using three nominal levels, i.e., $\alpha=0.01,$ $\alpha=0.05$ and $\alpha=0.1.$ All bootstrap calculations are based on $B=1000$ bootstrap replicates and $R=1000$ model repetitions. To examine the empirical size and power behavior of the TBB-based test, the curves in the two samples were generated according to model~\eqref{model:testing} and with the errors $\varepsilon_{i,t}$ following  model~\eqref{model:FAR1}, for $i\in\{1,2\}$, with mean functions given by $\mu_1(t) = 0$ and $\mu_2(t)=\gamma t(1-t)$ for the first and for the second population, respectively; see also Horv\'{a}th {\em et al.} (2013). The results obtained are shown in Table~\ref{tab:power} for a range of values of $\gamma.$ Notice that $\gamma=0$ corresponds  to the null hypothesis.

As it is evident from this table, the TBB-based test statistic $U_M^+$ has a good size behavior even in the case of $n_1=n_2=100$ observations while for $n_1=n_2=200$ observations the sizes of the TBB-based test are quite close to the nominal sizes for a range of block length values.
\begin{table}[htbp]
  \centering
    \begin{tabular}{|c|ccccc|cccc|}
    \hline
          &        \multicolumn{4}{|c}{$n_1=n_2=100$} &       & \multicolumn{4}{c|}{$n_1=n_2=200$} \\
    \multicolumn{1}{|c}{$\gamma$} & \multicolumn{1}{|c}{b} & \multicolumn{1}{l}{$\alpha=0.01$} & \multicolumn{1}{l}{$\alpha=0.05$} & \multicolumn{1}{l}{$\alpha=0.1$} &       &  \multicolumn{1}{c}{b} & \multicolumn{1}{l}{$\alpha=0.01$} & \multicolumn{1}{l}{$\alpha=0.05$} & \multicolumn{1}{l|}{$\alpha=0.1$} \\
\hline  0 & 4     & 0.026 & 0.077 & 0.142 &       &   6   &  0.013 & 0.057 & 0.113 \\
          & 6     & 0.015 & 0.061 & 0.112 &       &   8   &  0.010 & 0.052 & 0.115 \\
          & 8     & 0.015 & 0.071 & 0.128 &       &  10   &  0.013 & 0.066 & 0.106 \\
          & $b^*$ & 0.027 & 0.074 & 0.143 &       & $b^*$ &  0.013 & 0.057 & 0.113 \\
 &&&&&&&&&\\
    0.2 
          & 4     & 0.048 & 0.135 & 0.206 &       &   6   & 0.058 & 0.160 & 0.237 \\
          & 6     & 0.045 & 0.126 & 0.206 &       &   8   & 0.065 & 0.158 & 0.253 \\
          & 8     & 0.034 & 0.118 & 0.185 &       &  10   & 0.070 & 0.162 & 0.247 \\
          &$b^*$  & 0.042 & 0.116 & 0.178 &       & $b^*$ & 0.058 & 0.160 & 0.237 \\
 &&&&&&&&&\\
    0.5   
          & 4     & 0.225 & 0.418 & 0.544 &       &   6    & 0.408 & 0.615 & 0.715 \\
          & 6     & 0.200 & 0.374 & 0.499 &       &   8    & 0.411 & 0.632 & 0.759 \\
          & 8     & 0.184 & 0.356 & 0.490 &       &  10    & 0.425 & 0.645 & 0.749 \\
          &$b^*$  & 0.218 & 0.424 & 0.532 &       & $b^*$  & 0.408 & 0.615 & 0.715 \\
&&&&&&&&&\\
    0.8   
          & 4     & 0.584 & 0.772 & 0.853 &        & 6    & 0.864 & 0.966 & 0.980 \\
          & 6     & 0.543 & 0.763 & 0.841 &        & 8    & 0.865 & 0.948 & 0.975 \\
          & 8     & 0.529 & 0.739 & 0.831 &        & 10   & 0.843 & 0.948 & 0.976 \\
          &$b^*$  & 0.557 & 0.752 & 0.825 &       &$b^*$ & 0.864 & 0.966 & 0.980 \\
&&&&&&&&&\\
    1     
          & 4     & 0.779 & 0.898 & 0.945 &        &  6    & 0.972 & 0.995 & 0.998 \\
          & 6     & 0.746 & 0.891 & 0.941 &        &  8    & 0.975 & 0.994 & 0.999 \\
          & 8     & 0.755 & 0.898 & 0.943 &        &  10   & 0.969 & 0.994 & 0.998 \\
          &$b^*$  & 0.769 & 0.901 & 0.945 &        & $b^*$ & 0.972 & 0.995 & 0.998 \\
\hline
    \end{tabular}
\caption{Empirical size and power of the test based on TBB critical values and FAR(1) errors.}
   \label{tab:power}
\end{table}
It seems that the choice of the block size has a moderate effect on the power of the test.  Furthermore, the power of the TBB-based test increases as the deviations from the null become larger (i.e., larger values of $\gamma$) and/or as the sample size increases. Finally, using the suggested simple method to choose the block size $b,$  the corresponding test has good size and power behavior in all cases.

\subsection{\sc A real-life data example}
\label{sec:simulations2}
\begin{figure}
\begin{subfigure}{0.5\textwidth}
  \centering
  \includegraphics[trim = 10mm 15mm 10mm 15mm,clip=true,width=0.85\textwidth]{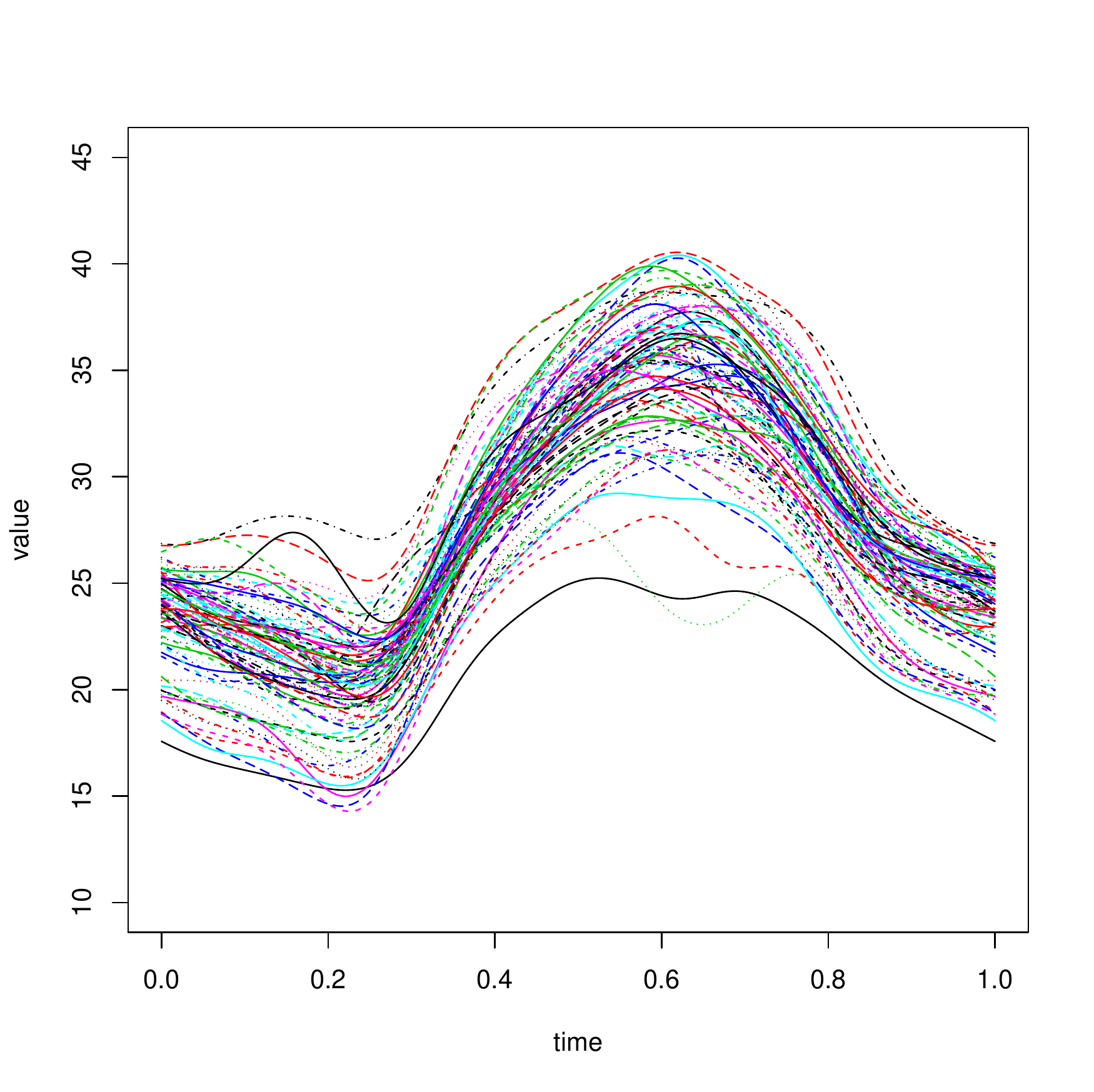}
\end{subfigure}%
\begin{subfigure}{0.5\textwidth}
  \centering
 \includegraphics[trim = 10mm 15mm 10mm 15mm,clip=true,width=0.85\textwidth]{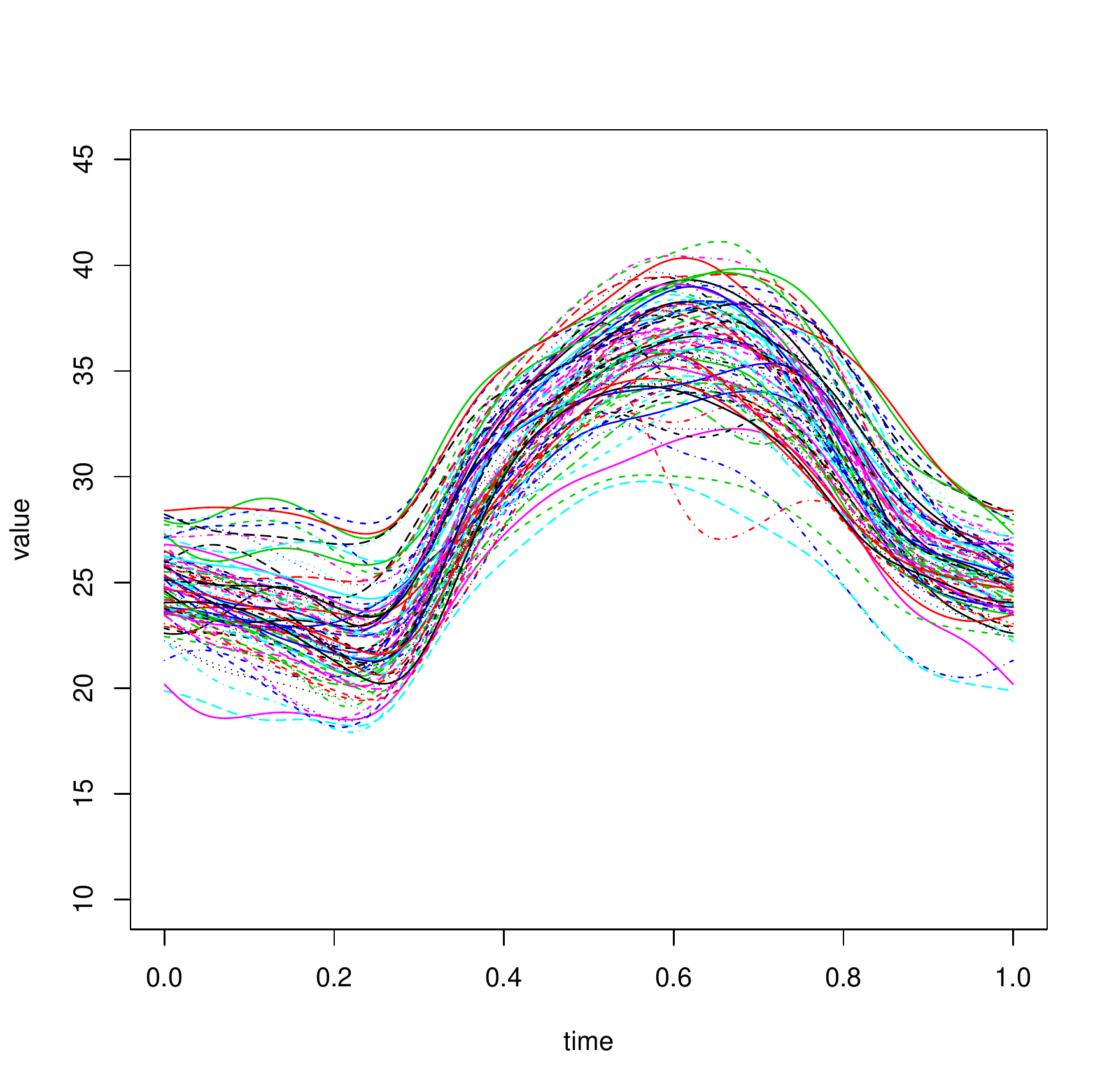}
\end{subfigure}
\caption{Temperature curves: summer 2005 (left  panel) and summer 2009 (right panel).}
\label{fig:summer}
\end{figure}

We apply the TBB-based testing procedure to a data set consisting of the summer temperature measurements recorded in Nicosia, Cyprus, for the years $2005$ and $2009.$ Our aim is to test whether there is a significant increase in the mean summer temperatures in $2009.$
The data consists of two samples of curves  $\{X_{i,t}(\tau),i=1,2,t=1,2,\ldots, 92\}$, where, $X_{i,t}$ represents the temperature of day $t$ of the summer $2005$ for $i=1$ and of the summer $2009$ for $i=2.$ More precisely, $X_{i,1}$ represents the temperature of the 1st of June and $X_{i,92}$ the temperature of the 31st of August. The temperature recordings have been taken in $15$ minutes intervals, i.e., there are $96$ temperature measurements for each day. These measurements have been transformed into functional objects using the Fourier basis with 21 basis functions. All curves are rescaled in order to be defined in the interval $\mathcal{I}.$ Figure~\ref{fig:summer} shows the temperatures curves of the summer of 2005 and of 2009. 

Since we are interested in checking whether there is an increase in the summer temperature in the year 2009 compared to 2005, the hypothesis of interest is $H_0:\mu_1(\tau) =\mu_2(\tau)$ versus $H_1:\mu_1(\tau)<\mu_2(\tau)$, for all $\tau \in \cal{I}$. The $p$-values of the TBB-based test using the test statistic $\widetilde{U}_M$ are:
0.001 (for $b=4$), 0.003 (for $b=6$),  0.004 (for $b=8$) and 0.002 (for $b=b^*$).
These $p$-values have been obtained using $B=1000$ bootstrap replicates.
As it is evident from these results, the $p$-values of the test statistic $\widetilde{U}_M$ are quite small leading to the rejection of $H_0$ for all commonly used $\alpha$-levels.

\section{\sc Appendix : Proofs}
To prove Theorem~\ref{thm:CLT} and Theorem \ref{thm:CLTtbb}, we first establish Lemma 5.1 and  Lemma 5.2. The proofs of Lemma 5.2 and Theorem \ref{thm:CLTtbb} are given in Section 1 of the supplementary material.
Note also that, throughout the proofs, we 
 use the fact that, by stationarity, $\E\|X_{i,m}-X_i\|=\E\|X_{0,m}-X_0\|$ and $\E\|X_{i,m}\|=\E\|X_i\|=\E\|X_0\|$ for all $i\in \mathbb{Z}$.
\begin{lemma}
Let $g_b$ be a non-negative, continuous and bounded function defined on $\mathbb{R}$, satisfying $g_b(0)=1$, $g_b(u)=g_b(-u)$, $g_b(u)\leq 1$ for all $u$, $g_b(u)=0,$ if $|u|>c,$ for some $c>0.$ Suppose that $(X_t, t \in \mathbb{Z})$ satisfies Assumption~\ref{as:L2} and $b=b(n)$ is a sequence of integers such that $b^{-1}+bn^{-1/2}=o(1)$ as $n\to\infty.$ Assume further that, for any fixed $u$, $g_b(u)\to 1$ as $n\to\infty.$ 
Then, as $n\to\infty$,
$$
\sum_{h=-b+1}^{b-1}g_b(h)\hat{\gamma}_h\overset{P}\to\sum_{i=-\infty}^{\infty}\E(\langle X_0,X_i\rangle),
$$
where
$\hat{\gamma}_h=\frac{1}{n}\sum_{i=1}^{n-|h|}\langle X_i,X_{i+|h|}\rangle$ for $-b+1\leq h\leq b-1$. \label{lemma:sigklisi}
\end{lemma}

\noindent\textbf{\textit{Proof.}}
First, note by the independence of $X_0$ and $X_{i,i}$, that
$
\sum_{i=1}^{\infty}|\E\langle X_0,X_i\rangle|=\sum_{i=1}^{\infty}|\E\langle X_0,X_i-X_{i,i}\rangle|\leq\sum_{i=1}^{\infty}(\E\|X_0\|^2)^{1/2}(\E\|X_0-X_{0,i}\|^2)^{1/2},
$
which implies by \eqref{orismosL2b} that $\sum_{i=-\infty}^{\infty}|\E\langle X_0,X_i\rangle|<\infty.$ Since $n^{-1}\sum_{i=1}^n\langle X_i,X_i\rangle-\E\langle X_0,X_0\rangle=o_P(1)$ as $n\to\infty$, it suffices to show that, as $n\to\infty$,
\begin{equation}
\label{SigklisiLemma}
\sum_{h=1}^{b-1}g_b(h)
\dfrac{1}{n}\sum_{i=1}^{n-h}\langle X_i,X_{i+h}\rangle-\sum_{i\geq1}\E\langle X_0,X_i\rangle=o_P(1).
\end{equation}
Let
$c_{\infty}^+= \sum_{i\geq1}\E[\langle X_0,X_i\rangle]$, $c_m^+ =\sum_{i\geq1}\E[\langle X_{0,m},X_{i,m}\rangle]$ and $ \hat{\gamma}_h^{(m)}=\dfrac{1}{n}\sum_{i=1}^{n-h}\langle X_{i,m},X_{i+h,m}\rangle.$ Because
\begin{equation}
\left|\sum_{h=1}^{b-1}g_b(h)\hat{\gamma}_h-c^+_\infty\right|\leq
|c_m^+-c^+_\infty|+\left|\sum_{h=1}^{b-1}g_b(h)\hat{\gamma}_h^{(m)}-c^+_m\right|
+\left|\sum_{h=1}^{b-1}g_b(h)\hat{\gamma}_h-\sum_{h=1}^{b-1}g_b(h)\hat{\gamma}_h^{(m)}\right|,
\label{eq:athrismata}
\end{equation}
assertion~\eqref{SigklisiLemma} is proved by showing that there exists $m_0\in\mathbb{N}$ such that all three terms on the right hand side of~\eqref{eq:athrismata} can be made arbitrarily small in probability as $n\to\infty$ for all $m \geq m_0$.\par
For the first term, we use the bound
\begin{align}
\left|\sum_{i\geq 1}\E\left[\langle X_{0,m},X_{i,m} \rangle -\langle X_0,X_i\rangle\right]\right|\leq&\left|\sum_{i=1}^m\E\left[\langle X_{0,m},X_{i,m} \rangle -\langle X_0,X_i\rangle\right]\right| \nonumber\\&\qquad+\left|\sum_{i=m+1}^\infty\E\left[\langle X_{0,m},X_{i,m} \rangle -\langle X_0,X_i\rangle\right]\right|,\label{eq:boundoffirstterm}
\end{align}
and handle the first term on the right hand side of~\eqref{eq:boundoffirstterm} using 
$
\langle X_{0,m},X_{i,m} \rangle -\langle X_0,X_i\rangle = \langle X_{0,m}- X_0,X_{i,m}\rangle+\langle X_0,X_{i,m}-X_i\rangle.
$
Cauchy-€"Schwarz's inequality and Assumption~\ref{as:L2} yields that for every $\epsilon_1 >0,$ $\exists\,m_1\in\mathbb{N}$ such that
\begin{align*}
\left|\sum_{i=1}^m\E\left[\langle X_{0,m},X_{i,m} \rangle -\langle X_0,X_i\rangle\right]\right|
& \leq 2\sum_{i=1}^{m}\left(\E\|X_{0,m}-X_0\|^2\E\|X_0\|^2\right)^{1/2}
\\& \leq 2(\E\|X_0\|^2)^{1/2} \left(m\left[\E\|X_{0,m}-X_0\|^2\right]^{1/2}\right)
<\epsilon_1
\end{align*}
for all $m\geq m_1$.
For the second term of the right hand side of~\eqref{eq:boundoffirstterm}, we get, using $\langle X_0,X_i \rangle = \langle X_{i,i},X_0 \rangle +\langle X_0,X_i-X_{i,i}\rangle,$ the fact that $X_0$ and $X_{i,i}$ as well as $X_{0,m}$ and $X_{i,m}$ are independent for $i\geq m+1$ and Lemma $2.1$ of Horv\'{a}th \& Kokoszka (2012), that, for any $\epsilon_2>0$, there exists $m_2 \in \mathbb{N}$ such that
\begin{align*}
\left|\sum_{i=m+1}^\infty\E\left[\langle X_{0,m},X_{i,m} \rangle -\langle X_0,X_i\rangle\right]\right|
&\leq \left|\sum_{i=m+1}^\infty\E\left[\langle X_{i,i},X_0\rangle\right]\right|+\left|\sum_{i=m+1}^\infty\E\left[\langle X_0,X_i-X_{i,i}\rangle\right]\right|\\& \leq \sum_{i=m+1}^\infty\left(\E\|X_0\|^2\E\|X_i-X_{i,i}\|^2\right)^{1/2}\\& =\left(\E\|X_0\|^2\right)^{1/2}\sum_{i=m+1}^\infty \left(\E\|X_0-X_{0,i}\|^2\right)^{1/2} < \epsilon_2
\end{align*}
for all $m \geq m_2$ because of \eqref{orismosL2b}. For the second term of \eqref{eq:athrismata}, first note that, for every fixed $m\geq1$ and for any fixed $h$, we have that
$
\left|\hat{\gamma}_h^{(m)}-\E[\langle X_{0,m},X_{h,m}\rangle]\right|=o_p(1).
$
Furthermore, since $\{X_{n,m},\,n \in \mathbb{Z}\}$ is an $m$-dependent sequence,
$
c_m^+=\sum_{i=1}^m\E[\langle X_{0,m},X_{i,m}\rangle].
$
Hence, the second term of the right hand side of~\eqref{eq:athrismata} is $o_p(1),$ if we show that
$|\sum_{h=m+1}^{b-1}g_b(h)\hat{\gamma}_h^{(m)}|=o_p(1)$.
We have 
$
\E\left[\sum_{h=m+1}^{b-1}g_b(h)\hat{\gamma}_h^{(m)}\right]^2= n^{-2}\sum_{h_1=m+1}^{b-1}\sum_{h_2=m+1}^{b-1}\sum_{i_1=1}^{n-h_1}\sum_{i_2=1}^{n-h_2} g_b(h_1)g_b(h_2)\E(\langle X_{i_1,m},X_{i_1+h_1,m}\rangle\langle X_{i_2,m},X_{i_2+h_2,m}\rangle).
$
Since the sequence $\{X_{i,m},\,i\in\mathbb{Z}\}$ is $m$-dependent, $X_{i,m}$ and $X_{i+h,m}$ are independent for $h\geq m+1,$ that is, using Lemma $2.1$ of Horv\'{a}th \& Kokoszka (2012) we have that $\E\langle X_{i,m},X_{i+h,m}\rangle=0$ for the same $h$. Hence, the number of non-vanishing terms $\E[\langle X_{i_1,m},X_{i_1+h_1,m}\rangle\langle X_{i_2,m},X_{i_2+h_2,m}\rangle]$ in the last equation above is of order $O(nb)$ and, consequently,
$
\E\left[\sum_{h=m+1}^{b-1}g_b(h)\hat{\gamma}_h^{(m)}\right]^2=O(b/n)=o(1)
$
from which the desired convergence follows by Markov's inequality. 
For the third term in~\eqref{eq:athrismata}, we show that, for $m \geq m_0$,
\begin{equation}
\label{eq:1osoros}
\limsup\limits_{n\rightarrow\infty}P\left(\left|\sum_{h=1}^{b-1} g_b(h)\left(\hat{\gamma}_h-\hat{\gamma}_h^{(m)}\right)\right|>\delta\right)=0,
\end{equation}
for all $\delta >0.$ From this, it suffices to show that, for $m \geq m_0$,
\begin{equation}
\label{eq:1osorosb}
\E\left|\sum_{h=1}^{b-1}g_b(h)(\hat{\gamma}_h-\hat{\gamma}_h^{(m)})\right|=o(1).
\end{equation}
Now, by the definitions of $\hat{\gamma}_h$ and $\hat{\gamma}_h^{(m)}$, we have
\begin{align}
\E\left|\sum_{h=1}^{b-1}g_b(h)\left(\hat{\gamma}_h-\hat{\gamma}_h^{(m)}\right)\right|& \leq \E\left|\dfrac{1}{n}\sum_{h=1}^m g_b(h)\sum_{i=1}^{n-h}\left(\langle X_i,X_{i+h}\rangle-\langle X_{i,m},X_{i+h,m}\rangle\right)\right|\nonumber\\& \quad
+\E\left|\dfrac{1}{n}\sum_{h=m+1}^{b-1} g_b(h) \sum_{i=1}^{n-h}\left(\langle X_i,X_{i+h}\rangle-\langle X_{i,m},X_{i+h,m}\rangle\right)\right|.\label{eq:lemmateliki}
\end{align}
For the first term of the right hand side of the above inequality, we use
$\langle X_i,X_{i+h}\rangle-\langle X_{i,m},X_{i+h,m}\rangle=\langle X_i-X_{i,m},X_{i+h}\rangle+\langle X_{i+h}-X_{i+h,m},X_{i,m}\rangle$,
and we get, by 
to get, by Cauchy-Schwarz's  inequality and simple algebra, that,
\begin{align*}
\E&\left|\dfrac{1}{n}\sum_{h=1}^m g_b(h) \sum_{i=1}^{n-h}(\langle X_i,X_{i+h}\rangle-\langle X_{i,m},X_{i+h,m}\rangle)\right|\\
&\quad
\leq m[(\E\| X_0-X_{0,m}\|^2\E\|X_0\|^2)^{1/2}+(\E\| X_0-X_{0,m}\|^2\E\|X_{0,m}\|^2)^{1/2}].
\end{align*}
Assumption~\ref{as:L2} implies then that, for every $\epsilon_3>0,$ there exists $m_3\in\mathbb{N}$ such that, for every $m\geq m_3$, the last quantity above is bounded by $\epsilon_3.$
For the second term on the right hand side of~\eqref{eq:lemmateliki}, we use the bound
\begin{equation}
\E\left|\dfrac{1}{n}\sum_{h=m+1}^{b-1}g_b(h) \sum_{i=1}^{n-h}\langle X_i,X_{i+h}\rangle\right|+\E\left|\dfrac{1}{n}\sum_{h=m+1}^{b-1}g_b(h) \sum_{i=1}^{n-h}\langle X_{i,m},X_{i+h,m}\rangle\right|.\label{eq:provlima1}
\end{equation}
Note  that the second summand of~\eqref{eq:provlima1} is $o(1)$, while for the first term we use 
$
\langle X_i,X_{i+h}\rangle=\langle X_i,X_{i+h,h}\rangle+\langle X_i,X_{i+h}-X_{i+h,h}\rangle
$
to get the bound
\begin{align}
\E\left|\dfrac{1}{n}\sum_{h=m+1}^{b-1}g_b(h)\sum_{i=1}^{n-h}\langle X_i,X_{i+h,h}\rangle\right|+\E\left|\dfrac{1}{n}\sum_{h=m+1}^{b-1}g_b(h) \sum_{i=1}^{n-h}\langle X_i,X_{i+h}-X_{i+h,h}\rangle\right|.\label{eq:provlima2}
\end{align}
For the last term of expression~\eqref{eq:provlima2}, we get, using~\eqref{orismosL2b}, that for every $\epsilon_4>0,$ there exists $m_4 \in \mathbb{N}$ such that
\begin{align*}
\dfrac{1}{n}\sum_{h=m+1}^{b-1} \sum_{i=1}^{n-h}\E\left|\langle X_i,X_{i+h}-X_{i+h,h}\rangle\right|&\leq
\sum_{h=m+1}^{b-1}\E\left|\langle X_0,X_{h}-X_{h,h}\rangle\right|\\& \leq(\E\|X_0\|^2)^{1/2}\sum_{h=m+1}^{\infty}(\E\|X_{0}-X_{0,h}\|^2)^{1/2} < \epsilon_4
\end{align*}
for all $m\geq m_4$.  Consider next the first term of~\eqref{eq:provlima2}. Because
$
\langle X_i,X_{i+h,h}\rangle=\langle X_i-X_{i,h},X_{i+h,h}\rangle+\langle X_{i,h},X_{i+h,h}\rangle,
$
we get for this term the bound
\begin{align}
\E\left|\dfrac{1}{n}\sum_{h=m+1}^{b-1}g_b(h)\sum_{i=1}^{n-h}\langle X_i-X_{i,h},X_{i+h,h}\rangle\right|+\E\left|\dfrac{1}{n}\sum_{h=m+1}^{b-1}g_b(h) \sum_{i=1}^{n-h}\langle X_{i,h},X_{i+h,h}\rangle\right|.\label{eq:provlima3}
\end{align}
The first term above is bounded by
$$
\E\left|\dfrac{1}{n}\sum_{h=m+1}^{b-1}g_b(h)\sum_{i=1}^{n-h}\langle X_i-X_{i,h},X_{i+h,h}\rangle\right|\leq(\E\|X_0\|^2)^{1/2}\sum_{h=m+1}^{\infty}(\E\|X_{0}-X_{0,h}\|^2)^{1/2}.
$$
Thus,  and by \eqref{orismosL2b}, for every $\epsilon_5>0,$ there exists $m_5 \in \mathbb{N}$ such that, for every $m\geq m_5$, this term is bounded by $\epsilon_5.$ For the last term of~\eqref{eq:provlima3}, note that $\{\langle X_{i,h},X_{i+h,h}\rangle,\,i\in\mathbb{Z}\}$ is an $2h$-dependent stationary process, and since  $X_i$ and $X_{i+h,h}$ are independent, i.e., $\E\langle X_i,X_{i+h,h}\rangle=0$ for all $i \in \mathbb{Z}$, $\{\langle X_{i,h},X_{i+h,h}\rangle,\,i\in\mathbb{Z}\}$ is then a mean zero $2h$-dependent stationary process which implies that $n^{-1/2}\sum_{i=1}^n\langle X_{i,h},X_{i+h,h}\rangle=O_P(1).$ Using Portmanteau's theorem, and since the function $f(x)=|x|$ is Lipschitz, we get that $\E\left|n^{-1/2}\sum_{i=1}^n\langle X_{i,h},X_{i+h,h}\rangle\right|=O(1).$ Therefore,
$$
\E\left|\dfrac{1}{n}\sum_{h=m+1}^{b-1}g_b(h)\sum_{i=1}^{n-h}\langle X_{i,h},X_{i+h,h}\rangle\right| \leq\dfrac{1}{\sqrt{n}}\sum_{h=m+1}^{b-1}\E\left|\dfrac{1}{\sqrt{n}}\sum_{i=1}^n\langle X_{i,h},X_{i+h,h}\rangle\right| =O(b/\sqrt{n})=o(1),
$$
which concludes the proof of the lemma by choosing $m_0=\max\{m_1,m_2,m_3,m_4,m_5\}$.

\begin{lemma}\label{lemmasTBB}$\,$
Suppose that $(Y_t, t \in \mathbb{Z})$ satisfies Assumption~\ref{as:L2} and that $b=b(n)$ is a sequence of integers satisfying $b^{-1}+bn^{-1/2}=o(1)$ as $n\to\infty.$ Let $w_n(\cdot),\,i=1,2,\ldots,$ be a sequence of data-tappering windows satisfying Assumption~\ref{as:giawtbb}. Then, as $n\to\infty$,
\begin{enumerate}[label=(\roman*),ref=\ref{lemmasTBB}\,(\roman*)]
\item $$
\sum_{|h|<b}\left(\dfrac{\mathcal{W}_{|h|}}{\|w_b\|^2_2}\right)E[\langle Y_0,y \rangle\langle Y_h,y \rangle]\to \sum_{i=-\infty}^{\infty} E[\langle Y_0,y \rangle\langle Y_i,y \rangle]\quad\text{for every}\;\; y\in L^2,
    $$ \label{lemma:sigklisi2TBB}
\item $$\iint\{\tilde{c}_n(u,v)-c(u,v)\}^2\mathrm{d}u\mathrm{d}v=o_P(1),
$$\label{lemma:sigklisiTBB}
\end{enumerate}
where $c(u,v)=\sum_{i=-\infty}^{\infty}\E[Y_0(u)Y_i(v)]$, $\mathcal{W}_h=\sum_{i=1}^{b-h}w_b(i)w_b(i+h),\,h=0,1,\ldots,b-1$
 and
$$\tilde{c}_n(u,v)=\dfrac{1}{n}\sum_{i=1}^{n} Y_i(u)Y_i(v)+
\sum_{h=1}^{b-1}\dfrac{\mathcal{W}_h}{\|w_b\|_2^2}
\dfrac{1}{n}\sum_{i=1}^{n-h}[Y_i(u)Y_{i+h}(v)+Y_{i+h}(u)Y_i(v)].$$
\end{lemma}

\noindent\textbf{\textit{Proof of Theorem~\ref{thm:CLT}.}}
By the triangle inequality and Theorem~$1$ of Horv\'{a}th {\em et al.} (2013), the assertion of the theorem is established if we show that, as $n\to\infty,$
\begin{equation}
\label{expr:CLT}
\sqrt{n}(\overline{X}^*_n-\E^*(\overline{X}^*_n)) \Rightarrow  \Gamma, \;\; \text{in probability},
\end{equation}
where $\Gamma$ is a Gaussian process in $L^2$ with mean $0$ and covariance operator $C$ with kernel $c(u,v)=\E(\Gamma(u)\Gamma(v))$ given for any $u,v\in[0,1]^2$ by
$$
c(u,v)=\E[X_0(u)X_0(v)]+\sum_{i\geq1}\E[X_0(u)X_i(v)]+\sum_{i\geq1}\E[X_0(v)X_i(u)].
$$
Using the notation $S_n^*=\sqrt{n}(\overline{X}^*_n-\E^*(\overline{X}^*_n)),$ it follows from Proposition $7.4.2$ of Laha and Rohatgi (1979) that, to prove~\eqref{expr:CLT}, it suffices to prove that,
\begin{enumerate}[label=(L\arabic*)]
  \item $\langle S_n^*,y\rangle\overset{d}\to N(0,\sigma^2(y))$ for every $y\in L^2$ where $\sigma^2(y)=\langle C(y),y\rangle,$ and that
  \item the sequence $\{S_n^*,n\in\mathbb{N}\}$ is tight.
\end{enumerate}
Consider $(L1).$ To establish the desired weak convergence, we prove that, as $n\to\infty$,
\begin{equation}
\label{eq:L1i}
\Var^*\left(\langle S_n^*,y\rangle\right)\overset{P}\to\sigma^2(y)
\end{equation}
and that
\begin{equation}
\label{eq:L1ii}
\dfrac{\langle S_n^*,y\rangle}{\sqrt{\Var^*(\langle S_n^*,y\rangle)}}\overset{d}\to N(0,1).
\end{equation}
Consider \eqref{eq:L1i} and notice that
$
S_n^*=\frac{1}{\sqrt{k}}\sum_{i=1}^{k}\left[U_i^*-\E^*(U_i^*)\right],
$
where $U_i^*=b^{-1/2}(X^*_{(i-1)b+1}+X^*_{(i-1)b+2}+\ldots+X^*_{ib}),\,i=1,2,\ldots,k.$
Due to the block bootstrap resampling scheme, the random variables $U_i^*,\,i=1,2,\ldots,k$ are i.i.d. Thus, using $\langle S_n^*,y \rangle=k^{-1/2}\sum_{i=1}^k [W_i^*-\E^*(W_i^*)]$, where $W_i^*=\langle U_i^*,y\rangle,\,i=1,2,\ldots,k$, we have
\begin{equation}
\Var^*\left(\langle S_n^*,y\rangle\right)=\E^*(W_1^*)^2-(\E^*(W_1^*))^2.\label{eq:diasporaW1}
\end{equation}
Let $\mu^*=\E^*(W_1^*)$ and $U_i=b^{-1/2}(X_i+X_{i+1},\ldots+X_{i+b-1}),\,i=1,2,\ldots,N.$ We then have that
\begin{equation}
\mu^*
=\dfrac{\sqrt{b}}{N}\left[\sum_{i=1}^{n}\langle X_i,y\rangle-\sum_{j=1}^{b-1}\left(1-\dfrac{j}{b}\right)[\langle X_j,y\rangle+\langle X_{n-j+1} ,y\rangle]\right].\label{eq:mesiW1}
\end{equation}
Therefore,
$\E^*(\mu^*)=0$.
Using
\begin{align*}
&\left[\sum_{i=1}^{n}\langle X_i,y\rangle-\sum_{j=1}^{b-1}\left(1-\dfrac{j}{b}\right)(\langle X_j,y\rangle+\langle X_{n-j+1} ,y\rangle)\right]^2\\& \quad=\sum_{i=1}^{n}\sum_{j=1}^{n}\langle X_i,y\rangle\langle X_j,y \rangle-2\sum_{i=1}^{n} \sum_{j=1}^{b-1}\left(1-\dfrac{j}{b}\right)\langle X_i,y\rangle[\langle X_j,y\rangle+\langle X_{n-j+1} ,y\rangle]\\&
\qquad+\sum_{i=1}^{b-1}\sum_{j=1}^{b-1}\left(1-\dfrac{i}{b}\right) \left(1-\dfrac{j}{b}\right)[\langle X_i,y\rangle+\langle X_{n-i+1} ,y\rangle][\langle X_j,y\rangle+\langle X_{n-j+1} ,y\rangle]
\end{align*}
we get,
\begin{equation}
\E(\mu^*)^2=
\dfrac{b}{N^2}\sum_{i=1}^{n}\sum_{j=1}^{n}\E[\langle X_i,y\rangle \langle X_j,y\rangle ]+O(b^2/n)=O(b^2/n)\label{eq:voithitikiLinder2},
\end{equation}
where the last equality follows since, by Kronecker's lemma,
\begin{align}
\dfrac{1}{n}\sum_{i=1}^{n}\sum_{j=1}^{n}\E[\langle X_i,y\rangle \langle X_j,y\rangle]&= \sum_{|h|<n}\left(1-\dfrac{|h|}{n}\right)\E[\langle X_0,y\rangle\langle X_h,y\rangle]\nonumber \\&
\to\iint c(u,v)y(u)y(v)\mathrm{d}u\mathrm{d}v\label{eq:sigklisidn}
\end{align}
as $n\to\infty$. Since $\E^*(\mu^*)=0$, \eqref{eq:voithitikiLinder2} implies that
$
\mu^*=O_P(b/\sqrt{n}). 
$

Consider next the first term of the right hand side of expression~\eqref{eq:diasporaW1}. For this term, we have
\begin{align}
\label{eq:anamenomeniW1^2}
\E^*(W_1^*)^2&=\dfrac{1}{N}\sum_{i=1}^{N}\langle U_i,y\rangle^2 \\ \nonumber
&=\dfrac{1}{N}\sum_{i=1}^{n}\langle X_i,y\rangle\langle X_i,y\rangle\\ \nonumber
&\;+
\sum_{h=1}^{b-1}\left(1-\dfrac{h}{b}\right)
\dfrac{1}{N}\sum_{i=1}^{n-h}[\langle X_i,y\rangle\langle X_{i+h},y\rangle+\langle X_{i+h},y\rangle\langle X_i,y\rangle]\\ \nonumber &
\;-\dfrac{1}{N}\sum_{s=1}^{b-1}\left(1-\dfrac{s}{b}\right)[\langle X_s,y \rangle \langle X_s,y \rangle +\langle X_{n-s+1},y \rangle \langle X_{n-s+1},y \rangle]\\ \nonumber &
\;-\dfrac{1}{N}\sum_{t=1}^{b-1}\sum_{j=1}^{b-t}\left(1-\dfrac{j+t}{b}\right)[\langle X_j,y\rangle\langle X_{j+t},y\rangle+\langle X_{n-j+1-t},y\rangle\langle X_{n-j+1},y\rangle\\ \nonumber &\hspace{125pt}
+\langle X_{j+t},y\rangle\langle X_{j},y\rangle+\langle X_{n-j+1},y\rangle\langle  X_{n-j+1-t},y\rangle].
\end{align}
Thus,
\begin{align*}
\E^*(W_1^*)^2&=
\dfrac{1}{N}\sum_{i=1}^{n}\langle X_i,y\rangle\langle X_i,y\rangle\\&\quad+
\sum_{h=1}^{b-1}\left(1-\dfrac{h}{b}\right)
\dfrac{1}{N}\sum_{i=1}^{n-h}[\langle X_i,y\rangle\langle X_{i+h},y\rangle+\langle X_{i+h},y\rangle\langle X_i,y\rangle]\\&\qquad+O_P(b/n)+O_P(b^2/n),
\end{align*}
from which we get
\begin{equation}
\Var^*(W_1^*)=\iint c_N(u,v)y(u)y(v)\mathrm{d}u\mathrm{d}v+O_p(b^2/n),\label{eq:varW1}
\end{equation}
where
\begin{equation}
\label{eq:cNuv}
c_N(u,v)=\dfrac{1}{N}\sum_{i=1}^{n} X_i(u)X_i(v)+
\sum_{h=1}^{b-1}\left(1-\dfrac{h}{b}\right)
\dfrac{1}{N}\sum_{i=1}^{n-h}\Big[X_i(u)X_{i+h}(v)+X_{i+h}(u)X_i(v)\Big].
\end{equation}
By the ergodic theorem and equation $(A.2)$ of Horv\'{a}th {\em et al.} (2013), choosing the kernel $K$ in their notation to be the kernel $K(x)=(1-|x|)\mathds{1}_{[-1,1]}(x)$, where $\mathds{1}_{A}(x)$ denotes the indicator function of $A$, it follows that
\begin{equation}\label{eq:fromA2}
\iint [c_n(u,v)-c(u,v)]^2\mathrm{d}u\mathrm{d}v=o_P(1)
\end{equation}
as $n\to\infty$, where $c(u,v)=\sum_{i=-\infty}^{\infty}\E[X_0(u)X_i(v)]$ and
$c_n(u,v)=(N/n)c_N(u,v).$
Using Cauchy-Schwarz's inequality, we get that, as $n \to \infty$,
$$
\left|\iint (c_n(u,v)-c(u,v))y(u)y(v)\mathrm{d}u\mathrm{d}v\right|
\leq \left(\iint \{c_n(u,v)-c(u,v)\}^2\mathrm{d}u\mathrm{d}v\right)^{1/2}\|y\|^2
=o_P(1).
$$
That is,
$$
\iint c_n(u,v)y(u)y(v)\mathrm{d}u\mathrm{d}v\overset{P}\to\iint c(u,v)y(u)y(v)\mathrm{d}u\mathrm{d}v.
$$
Since $N/n\to 1$ as $n\to\infty$, we finally get from~\eqref{eq:varW1} that,
\begin{equation}
\Var^*\langle S^*_n,y \rangle=\Var^*(W^*_1)\overset{P}\to\iint c(u,v)y(u)y(v)\mathrm{d}u\mathrm{d}v=\sigma^2(y)\label{eq:sigklisiVarw1}.
\end{equation}

Consider next \eqref{eq:L1ii}. Observe that $W_i^*=\langle U_i^*,y\rangle,\,i=1,2,\ldots,k$ are i.i.d. random variables and, therefore, it suffices to show that Lindeberg's condition
\begin{equation}
\label{eq:Lindeberg}
\lim_{n\to\infty}\dfrac{1}{\tau_k^{*2}}\sum_{t=1}^{k}\E^*\left[(W_t^*-\mu^*)^2 \mathds{1}(|W_t^*-\mu^*|>\varepsilon \tau_k^*)\right]=0,\quad\text{for every}\,\,\varepsilon>0,
\end{equation}
is fulfilled, where
$\tau_k^{*2}=\sum_{t=1}^{k}\Var^*(W_t^*)=k\Var^*(W_1^*)$
and $\mu^*=\E^*(W_i^*).$ To establish~\eqref{eq:Lindeberg}, and because of~\eqref{eq:sigklisiVarw1}, it suffices to show that, for any $\delta>0$ and as $n\to\infty,$
\begin{equation}
\label{eq:sufficeLindeberg}
P\left(\dfrac{1}{k}\sum_{t=1}^{k}\E^*\left[(W_t^*-\mu^*)^2 \mathds{1}(|W_t^*-\mu^*|>\varepsilon \tau_k^*)\right]>\delta\right)\to0.
\end{equation}
Towards this, notice first that, for any two random variables $X$ and $Y$ and any $\eta>0$, it yields that
\begin{align}
\E[|X+Y|^2&\mathds{1}(|X+Y|>\eta)]\nonumber\\&\leq 
4\left[\E|X|^2\mathds{1}(|X|>\eta/2)+\E|Y|^2\mathds{1}(|Y|>\eta/2)\right];\label{eq:Lahiri}
\end{align}
see Lahiri $(2003)$, p. $56$. We then get by Markov's inequality that
\begin{align}
P&\left(\dfrac{1}{k}\sum_{t=1}^{k}\E^*\left[(W_t^*-\mu^*)^2 \mathds{1}(|W_t^*-\mu^*|> \varepsilon \tau_k^*)\right]>\delta\right) \nonumber \\&\quad
\leq\delta^{-1}\E\left\{\E^*\left[(W_1^*-\mu^*)^2 \mathds{1}(|W_1^*-\mu^*|> \varepsilon \tau_k^*)\right]\right\}\nonumber \\&\quad
= \delta^{-1}\E\left[\dfrac{1}{N}\sum_{i=1}^{N}(W_i-\mu^*)^2 \mathds{1}(|W_i-\mu^*|> \varepsilon \tau_k^*)\right]
\nonumber \\&\quad = \delta^{-1}\E\left[(W_1-\mu^*)^2 \mathds{1}(|W_1-\mu^*|> \varepsilon \tau_k^*)\right]\nonumber\\
&\quad\leq 4\delta^{-1}\left[\E W^2_1 \mathds{1}(|W_1|> \varepsilon \tau_k^*/2)+\E(\mu^*)^2\right],\label{eq:voithitikiLinder}
\end{align}
where $W_i=\langle U_i,y\rangle,\, i=1,2,\ldots,N.$ Furthermore, we have
$$
\E(W_1^2)=\E|\langle U_1,y \rangle|^2
=\sum_{|h|<b}\left(1-\dfrac{|h|}{b}\right)\E[\langle X_0,y \rangle\langle X_h,y \rangle] \to\iint c(u,v)y(u)y(v)\mathrm{d}u\mathrm{d}v,
$$
as $n \rightarrow \infty$.
Therefore, by the dominated convergence theorem,
$
\lim_{n\to\infty}\E W^2_1 \mathds{1}(|W_1|> \varepsilon \tau_k^*/2)=0,
$
Hence, using expression~\eqref{eq:voithitikiLinder2}, we conclude that~\eqref{eq:voithitikiLinder} converges to $0$ as $n\to\infty.$\par
To establish $(L2),$ it suffices, by Theorem $1.13$ of Prokhrov (1956) and Theorems $5.1$ and $5.2$ of Billingsley (1999), 
to prove that $\lim_{k \to \infty} \sup_{n\geq1}\sum_{j=k}^{\infty}\E|\langle S^*_n,e_j\rangle|^2 = 0,$ where $\{e_j,\,j\geq1\}$ is a complete orthonormal basis of $L^2$. Using
$
\E^*|\langle S^*_n,e_j\rangle|^2=\Var^*(\langle U_1^*,e_j\rangle)
$
and Lemma $14$ of Cerovecki and H\"{o}rmanm (2017), $(L2)$ is satisfied if the following five conditions are fulfilled.
\begin{enumerate}[label=(\alph*)]
  \item $\Var^*(\langle U_1^*,e_j\rangle)\geq0\quad\forall j,n;$
  \item $\lim_{n\to\infty}\Var^*(\langle U_1^*,e_j\rangle)=\Sigma_j,\quad\text{in probability};$
  \item $\sum_{j\geq 1}\Sigma_j<\infty;$
  \item $\lim_{n \to \infty}\sum_{j\geq 1}\Var^*(\langle U_1^*,e_j\rangle)= \sum_{j\geq 1}\Sigma_j,\quad\text{in probability};$
  \item $\sum_{j\geq1}\Var^*(\langle U_1^*,e_j\rangle)$ is bounded for all $n\geq1$, in probability.
\end{enumerate}
Note that, by letting $y=e_j$ in expression~\eqref{eq:sigklisiVarw1}, property (b) follows with $\Sigma_j=\iint c(u,v)e_j(u)e_j(v)\mathrm{d}u\mathrm{d}v.$ To prove (c), notice that, by Proposition $6$ of H\"{o}rmanm {\em et al.} (2015), and since the stochastic process $\{X_t,t\in\mathbb{Z}\}$ is $L^2$-$m$-approximable, the covariance operator $C$ with kernel $c(\cdot,\cdot)$ is trace class. Therefore,
\begin{equation}
\label{eq:conditionC}
\sum_{j\geq 1}\Sigma_j=\sum_{j\geq 1}\iint c(u,v)e_j(u)e_j(v)\mathrm{d}u\mathrm{d}v=\sum_{j\geq 1}\lambda_j<\infty,
\end{equation}
where $\lambda_j,\,j\geq1$ are the eigenvalues of $C.$

To establish (d), we get, using~\eqref{eq:mesiW1}, that
\begin{equation}
\Var^*(\langle U_1^*,e_j\rangle)= \dfrac{1}{N}\sum_{i=1}^{N}\langle U_i,e_j\rangle^2-\left(\dfrac{\sqrt{b}}{N}\left[\sum_{i=1}^{n}\langle X_i,e_j\rangle-\sum_{l=1}^{b-1}\left(1-\dfrac{l}{b}\right)[\langle X_l,e_j\rangle+\langle X_{n-l+1} ,e_j\rangle]\right]\right)^2. \label{eq:d1}
\end{equation}
By Parseval's identity, we have,
\begin{align*}
\sum_{j=1}^{\infty}\dfrac{1}{N}&\sum_{i=1}^{N}|\langle U_i,e_j\rangle|^2=\dfrac{1}{N}\sum_{i=1}^{N}\|U_i\|^2\\&
=\dfrac{1}{N}\sum_{i=1}^{n}\langle X_i,X_i\rangle+\sum_{h=1}^{b-1}\left(1-\dfrac{h}{b}\right)\dfrac{1}{N}\sum_{i=1}^{n-h}[\langle X_i,X_{i+h}\rangle+\langle X_{i+h},X_i\rangle]\\&\hspace{20pt}
-\dfrac{1}{N}\sum_{s=1}^{b-1}\left(1-\dfrac{s}{b}\right)[\langle X_s,X_s \rangle +\langle X_{n-s+1},X_{n-s+1}\rangle]\\&\hspace{20pt}
-\dfrac{1}{N}\sum_{t=1}^{b-1}\sum_{j=1}^{b-t}\left(1-\dfrac{t+j}{b}\right)
[\langle X_j,X_{j+t}\rangle+\langle X_{n-j+1-t},X_{n-j+1}\rangle\\&\hspace{145pt}+
\langle X_{j+t},X_{j}\rangle+\langle X_{n-j+1},X_{n-j+1-t}\rangle].
\end{align*}
Hence,
\begin{equation}
\sum_{j=1}^{\infty}\dfrac{1}{N}\sum_{i=1}^{N}|\langle U_i,e_j\rangle|^2=\dfrac{1}{N}\sum_{i=1}^{n}\langle X_i,X_i\rangle+ \sum_{h=1}^{b-1}\left(1-\dfrac{h}{b}\right)\dfrac{1}{N}\sum_{i=1}^{n-h}[\langle X_i,X_{i+h}\rangle+\langle X_{i+h},X_i\rangle]+O_P(b^2/n)\label{eq:d2}.
\end{equation}
Then, by letting $g_b(h)=\left(1-\frac{|h|}{b}\right)$ in Lemma~\ref{lemma:sigklisi}, we get that, as $n\to\infty$,
\begin{equation}
\label{eq:sigkisid}
\sum_{j=1}^{\infty}\dfrac{1}{N}\sum_{i=1}^{N}\langle U_i,e_j\rangle^2\overset{P}\to\sum_{i=-\infty}^{\infty}\E(\langle X_0,X_i\rangle).
\end{equation}
For the second term of equation~\eqref{eq:d1}, we show that,
\begin{equation}
\label{eq:sinthikigiad}
\sum_{j\geq 1}\left(\dfrac{\sqrt{b}}{N}\left[\sum_{i=1}^{n}\langle X_i,e_j\rangle-\sum_{l=1}^{b-1}\left(1-\dfrac{l}{b}\right)[\langle X_l,e_j\rangle+\langle X_{n-l+1} ,y\rangle]\right]\right)^2=o_P(1),
\end{equation}
as $n\to\infty.$ Using $\langle x,y \rangle=\sum_{j\geq 1}\langle x,e_j\rangle\langle y,e_j\rangle,$ we have
\begin{align*}
\sum_{j\geq 1}&\left[\sum_{i=1}^{n}\langle X_i,e_j\rangle-\sum_{l=1}^{b-1}\left(1-\dfrac{l}{b}\right)(\langle X_l,e_j\rangle+\langle X_{n-l+1} ,e_j\rangle)\right]^2\\
&=
\sum_{i=1}^{n}\sum_{l=1}^{n}\langle X_i,X_l \rangle-2\sum_{i=1}^{n} \sum_{l=1}^{b-1}\left(1-\dfrac{l}{b}\right)\left[\langle X_i,X_l\rangle+\langle X_i,X_{n-l+1}\rangle\right]\\&
\qquad+\sum_{i=1}^{b-1}\sum_{l=1}^{b-1}\left(1-\dfrac{i}{b}\right) \left(1-\dfrac{l}{b}\right)\left[\langle X_i,X_l\rangle+\langle X_{n-i+1} ,X_l\rangle+\langle X_i,X_{n-l+1}\rangle+\langle X_{n-i+1} ,X_{n-l+1}\rangle\right]\\&=
\sum_{i=1}^{n}\sum_{l=1}^{n}\langle X_i,X_l \rangle+O_P(nb)+O_P(b^2).
\end{align*}
Now note that, by the continuous mapping theorem and using Theorem~$1$ of Horv\'{a}th {\em et al.} (2013), we get
\begin{equation}
\label{eq:antilemma}
\dfrac{1}{n}\sum_{i=1}^{n}\sum_{l=1}^{n}\langle X_i,X_l\rangle=\langle \sqrt{n}\overline{X_n},\sqrt{n}\overline{X_n},\rangle=O_p(1).
\end{equation}
Therefore, 
$$
\dfrac{b}{N^2}\left[\sum_{i=1}^{n}\sum_{l=1}^{n}\langle X_i,X_l \rangle+O_P(nb)+O_P(b^2)\right]=O_P(b^2/n)=o_p(1),
$$
which establishes~\eqref{eq:sinthikigiad}.
Hence, from~\eqref{eq:d1},~\eqref{eq:sigkisid} and~\eqref{eq:sinthikigiad}, we conclude that
\begin{equation}
\label{eq:voithitikigiae}
\sum_{j\geq 1}\Var^*(\langle U_1^*,e_j\rangle)\to \sum_{i=-\infty}^{\infty}\E(\langle X_0,X_i\rangle),\quad\text{in probability.}
\end{equation}
Therefore, and by~\eqref{eq:conditionC}, property (d) is proved if we show that,
\begin{equation}
\label{eq:MercerV}
\sum_{j\geq 1}\lambda_j=\sum_{i=-\infty}^{\infty}\E(\langle X_0,X_i\rangle).
\end{equation}
Using Mercer's theorem, we have
\begin{equation}
\sum_{j\geq 1}\lambda_j=\int c(u,u) \mathrm{d}u=\int\sum_{i=-\infty}^\infty\E[X_0(u)X_i(u)]\mathrm{d}u 
=\sum_{i=-\infty}^\infty\E\int[X_0(u)X_i(u)]\mathrm{d}u = \sum_{i=-\infty}^{\infty}\E\langle X_0,X_i\rangle.\label{eq:Mercer}
\end{equation}
Notice that the above interchange of summation and integration is justified since, using Assumption~\ref{as:L2}, and the fact that $X_0$ and $X_{i,i}$ are independent for $i\geq 1$, 
we get
\begin{align*}
\sum_{i=-\infty}^{\infty}&\int\left|\E [X_0(u)X_i(u)]\right|\mathrm{d}u
\\&= \int\left|\E[X_0(u)X_0(u)]\right|\mathrm{d}u+2\sum_{i=1}^\infty\int\left|\E \{X_0(u)[X_i(u)-X_{i,i}(u)]\}\right|\mathrm{d}u
\\&\leq\int\E(X_0(u))^2\mathrm{d}u+2\sum_{i=1}^{\infty}\left\{\int\E[X_0(u)]^2\mathrm{d}u\right\}^{1/2} \left\{\int\E[X_i(u)-X_{i,i}(u)]^2\mathrm{d}u\right\}^{1/2}
\\&\leq
\E\|X_0\|^2+2\left(\E\|X_0\|^2\right)^{1/2} \sum_{i=1}^{\infty} \left(\E\| X_0-X_{0,i}\|^2\right)^{1/2}<\infty.
\end{align*}

To prove (e), notice first that, by~\eqref{eq:d1}, $\sum_{j=1}^{\infty}\Var^*(\langle U_1^*,e_j\rangle)\leq\sum_{j=1}^{\infty}(1/N)\sum_{i=1}^{N}|\langle U_i,e_j\rangle|^2$ and, therefore, using~\eqref{eq:d2}, for any given $n \geq 1$, $\sum_{j=1}^{\infty}\Var^*(\langle U_1^*,e_j\rangle)$ is bounded in probability. Furthermore, by~\eqref{eq:voithitikigiae}, the sequence $\{\sum_{j=1}^{\infty}\Var^*(\langle U_1^*,e_j\rangle),\,n\geq1\}$ converges in probability as $n\to\infty.$\par

Consider next assertion $(\rom{2})$ of the theorem. By the triangle inequality,  it suffices to prove that as $n\to\infty,$ $\|n\E^*(\overline{X}_n^*-\E^*(\overline{X}_n^*))\otimes(\overline{X}_n^*-\E^*(\overline{X}_n^*))-2\pi\mathcal{F}_0\|_{HS}=o_P(1).$ Now, recall that $U_i^*,\,i=1,2,\ldots,n$, are i.i.d., and note that
$$
n\E^*(\overline{X}_n^*-\E^*(\overline{X}_n^*))\otimes(\overline{X}_n^*-\E^*(\overline{X}_n^*))(y)(v)=\int\E^*\Big[[U_1^*(u)-\E^*(U^*_1(u))][U_1^*(v)-\E^*(U^*_1(v))]\Big]y(u)\mathrm{d}u, 
$$
i.e., $n\E^*(\overline{X}_n^*-\E^*(\overline{X}_n^*))\otimes(\overline{X}_n^*-\E^*(\overline{X}_n^*))$
is an integral operator with kernel
\begin{equation}
\label{eq:kerneld}
d(u,v)=\E^*[U_1^*(u)U_1^*(v)]-\E^*(U^*_1(u))\E^*(U^*_1(v)).
\end{equation}
Now,
\begin{align}
\E^*[U_1^*(u)U_1^*(v)]=&
\dfrac{1}{N}\sum_{i=1}^{n}X_i(u)X_i(v)+
\sum_{h=1}^{b-1}\left(1-\dfrac{h}{b}\right)
\dfrac{1}{N}\sum_{i=1}^{n-h}[X_i(u)X_{i+h}(v)+X_{i+h}(u)X_i(v)]\nonumber\\&
\;-\dfrac{1}{N}\sum_{s=1}^{b-1}\left(1-\dfrac{s}{b}\right)[X_s(u)X_s(v)+X_{n-s+1}(u)X_{n-s+1}(v)]\nonumber\\&
\;-\dfrac{1}{N}\sum_{t=1}^{b-1}\sum_{j=1}^{b-t}\left(1-\dfrac{j+t}{b}\right)[X_j(u)X_{j+t}(v)+X_{n-j+1-t}(u)X_{n-j+1-t}(v)\nonumber\\&\hspace{120pt}
+X_{j+t}(u) X_{j}(v)+X_{n-j+1}(u)X_{n-j+1-t}(v)]\label{eq:kerneld1}
\end{align}
and
\begin{equation}
\E^*(U^*_1(u))=\dfrac{\sqrt{b}}{N}\left[\sum_{i=1}^{n} X_i(u)-\sum_{j=1}^{b-1}\left(1-\dfrac{j}{b}\right)(X_j(u)+X_{n-j+1}(u))\right].\label{eq:kerneld2}
\end{equation}
Therefore, $d(u,v)=c_N(u,v)+R(u,v)$, where $R(u,v)$ is defined as the difference of $d(u,v)$ given in \eqref{eq:kerneld} and $c_N(u,v)$ given in \eqref{eq:cNuv}. Now, notice that
$
2\pi\mathcal{F}_0(y)(v)=\int\sum_{h=-\infty}^{\infty}\E[X_0(u)X_h(v)]y(u)\mathrm{d}u,
$
i.e., $2\pi\mathcal{F}_0$ is an integral operator with kernel $c(u,v)=\sum_{h=-\infty}^{\infty}\E[X_0(u)X_h(v)].$
Hence, 
\begin{align*}
\|n\E^*&(\overline{X}_n^*-\E^*(\overline{X}_n^*))\otimes(\overline{X}_n^*-\E^*(\overline{X}_n^*))-2\pi\mathcal{F}_0\|_{HS}\\&=
\iint[d(u,v)-c(u,v)]^2\mathrm{d}u\mathrm{d}v\leq 2\iint[c_N(u,v)-c(u,v)]^2\mathrm{d}u\mathrm{d}v+2\iint[R(u,v)]^2\mathrm{d}u\mathrm{d}v.
\end{align*}
Using~\eqref{eq:fromA2} it suffices to prove that $\iint[R(u,v)]^2\mathrm{d}u\mathrm{d}v=o_p(1).$ To prove this, recall the inequality $(\sum_{i=1}^{L}a_i)^2\leq L \sum_{i=1}^{L}a_i^2$, where $L$ is a positive integer, and notice that, using~\eqref{eq:antilemma},
\begin{align}
\dfrac{b^2}{N^4}\iint&\left(\sum_{i=1}^{n}\sum_{j=1}^{n}X_i(u)X_j(v)\right)^2\mathrm{d}u\mathrm{d}v\nonumber\\
&=\dfrac{b^2}{N^2}\dfrac{1}{N}\sum_{i_1=1}^{n}\sum_{i_2=1}^{n}\int X_{i_1}(u)X_{i_2}(u)\mathrm{d}u\dfrac{1}{N}\sum_{j_1=1}^{n}\sum_{j_2=1}^{n}\int X_{j_1}(v)X_{j_2}(v)\mathrm{d}v\nonumber\\&=\dfrac{b^2}{N^2}\left(\dfrac{1}{N}\sum_{i_1=1}^{n}\sum_{i_2=1}^{n}\langle X_{i_1},X_{i_2}\rangle\right)^2=O_P(b^2/N^2)=o_p(1)\label{eq:finalAs2a}.
\end{align}
Furthermore, 
\begin{align}
\iint\left[\dfrac{1}{N}\sum_{t=1}^{b-1}\sum_{j=1}^{b-t}\left(1-\dfrac{j+t}{b}\right)X_j(u)X_{j+t}(v)\right]^2\mathrm{d}u\mathrm{d}v &\leq\dfrac{1}{N^2}b^2\iint\sum_{t=1}^{b-1}\sum_{j=1}^{b-t}X^2_j(u)X^2_{j+t}(v)\mathrm{d}u\mathrm{d}v\nonumber\\&=O_P(b^4/N^2)=o_p(1),\label{eq:finalAs2b}
\end{align}
where all other terms appearing in $R(u,v)$ are handled similarly. This completes the proof of the theorem.

\noindent\textbf{\textit{Proof of Theorem~\ref{thm:consistent}.}}
Consider assertion $(\rom{1}).$ For $i=1,2$, let $\{e^*_{i,j},\,j=1,2,\ldots,n_i\}$ be the pseudo-observations generated by implementing the MBB procedure at $\{\varepsilon_{i,j},\,j=1,2,\ldots,n_i\}.$ Using Theorem~\ref{thm:CLT}, it follows that, conditionally on $\boldsymbol{\mathrm{X_M}}$, for $i=1,2$, and as $n_1, n_2\to\infty$,
$$
\frac{1}{\sqrt{n_i}}\sum_{j=1}^{n_i}(e^*_{i,j}-\E^*(e^*_{i,j})) \Rightarrow \Gamma_i, \;\; \text{in probability},
$$
where $\Gamma_i$ is a Gaussian random element with mean zero and covariance operator $C_i$ with kernel $c_i(\cdot,\cdot)$.
Now, recall from Step~$3$ of the MBB-based testing algorithm that, for $i=1,2$, the pseudo-observations $\varepsilon^*_{i,\xi+sb}(\tau),\,\xi=1,2,\ldots,b,\,s=0,1,\ldots,q_i,\,\,\tau\in\mathcal{I}$, are generated by first applying the MBB procedure to $\hat{\varepsilon}_{i,\xi+sb}(\tau),\,\xi=1,2,\ldots,b,\,s=0,1,\ldots,q_i,\,\,\tau\in\mathcal{I}$ and then $\overline{\varepsilon}_{i,\xi}(\tau)$ is subtracted.  Note further that $\varepsilon_{i,j}(\tau)=\hat{\varepsilon}_{i,j}(\tau)+\overline{X}_{i,n_i}-\mu_i(\tau).$
Thus, $e^*_{i,\xi+sb}(\tau)=\varepsilon^*_{i,\xi+sb}(\tau)+\overline{\varepsilon}_{i,\xi}(\tau)+\overline{X}_{i,n_i}(\tau)-\mu_i(\tau)$ and, using expression~\eqref{eq:mesialgorithmouMBB}, we get
$$
\dfrac{1}{\sqrt{n_i}}\sum_{j=1}^{n_i}(e^*_{i,j}-\E^*(e^*_{i,j}))=
\dfrac{1}{\sqrt{n_i}}\sum_{j=1}^{n_i}(X^*_{i,j}-\E^*(X^*_{i,j}))=
\dfrac{1}{\sqrt{n_i}}\sum_{j=1}^{n_i}(X^*_{i,j}-\overline{X}_M).
$$
Therefore, and conditionally on $\boldsymbol{\mathrm{X_M}}$, as $n_1,n_2\to\infty$,
$$
\left(\dfrac{1}{\sqrt{n_1}}\sum_{j=1}^{n_1}(X^*_{1,j}-\overline{X}_M), \dfrac{1}{\sqrt{n_2}}\sum_{j=1}^{n_2}(X^*_{2,j}-\overline{X}_M)\right) \Rightarrow (\Gamma_1,\Gamma_2), \;\; \text{in probability},
$$ 
where $\Gamma_1$ and $\Gamma_2$ are two independent Gaussian random elements with mean zero and covariance operator $C_1$ and $C_2$ with kernel $c_1(\cdot,\cdot)$ and $c_2(\cdot,\cdot)$, respectively. Since
$$
\sqrt{\dfrac{n_1n_2}{M}}(\overline{X}^*_{1,n_1}-\overline{X}^*_{2,n_2})=
\sqrt{\dfrac{n_2}{M}}\dfrac{1}{\sqrt{n_1}}\sum_{j=1}^{n_1}(X^*_{1,j}-\overline{X}_M)-\sqrt{\dfrac{n_1}{M}}\dfrac{1}{\sqrt{n_2}}\sum_{j=1}^{n_2}(X^*_{2,j}-\overline{X}_M),
$$
and because $n_1/M\to\theta$, we get that, as $n_1,n_2\to\infty$, 
$$
\sqrt{\frac{n_1n_2}{ M}} (\overline{X}^*_{1,n_1}-\overline{X}^*_{2,n_2}) \Rightarrow \Gamma, \;\; \text{in probability},
$$
where $\Gamma =\sqrt{1-\theta}\Gamma_1-\sqrt{\theta}\Gamma_2.$  
The proof of assertion $(\rom{2})$ follows along the same lines using Theorem~\ref{thm:CLTtbb}. This completes the proof of the theorem.

\section*{Acknowledgement}
We would like to thank the two referees for their useful comments and suggestions which have improved the presentation of the paper.

\section*{Supplementary Material}
The online supplement contains the proofs of Lemma 5.1 and Theorem 2.2 as well as some additional numerical results. You can download the file from:
\begin{center}
http://www.mas.ucy.ac.cy/$\sim$fanis/Papers/MBB-TBB-Supplement-Rev.pdf
\end{center}

\end{document}